\documentclass[11pt,a4paper]{article}

\usepackage{color}
\usepackage{graphics,epsfig}
\usepackage{amsfonts,amssymb,amsmath,amsthm}

\newtheorem{lem}{Lemma}[section]
\newtheorem{cor}[lem]{Corollary}
\newtheorem{thm}[lem]{Theorem}
\newtheorem{prop}[lem]{Proposition}
\theoremstyle{definition}
\newtheorem{defi}[lem]{Definition}
\theoremstyle{remark}
\newtheorem{rem}[lem]{Remark}

\numberwithin{equation}{section}

\newcommand{\surm}{\di{\frac{1}{m(x)}}}
\newcommand{\ep}{\varepsilon}
\newcommand{\ue}{u^\ep}
\newcommand{\U}{U_0}
\newcommand{\UU}{U_{0z}}
\newcommand{\UUU}{U_{0zz}}
\newcommand{\VV}{U_{1z}}
\newcommand{\VVV}{U_{1zz}}
\newcommand{\n}{\nabla }
\newcommand{\om}{\Omega}
\newcommand{\ombar}{\overline{\Omega}}
\newcommand{\edeux}{\displaystyle{\frac{1}{\ep ^2}}}
\newcommand{\R}{\mathbb{R}}

\newcommand{\vsp}{\vspace{8pt}}
\renewcommand{\div}{\operatorname{div}}
\newcommand{\di}{\displaystyle}
\newcommand{\EB}{e^{-\beta t/\ep ^2}}
\newcommand{\Pe}{(P^{\;\!\ep})}
\newcommand {\Q}{Q_T}
\newcommand{\emuth}{e^{\mu t/\ep ^2}}
\newcommand{\um}{u _\ep ^-}
\newcommand{\upl}{u _\ep ^+}
\newcommand{\dtip}{\widetilde d _\phi}
\newcommand{\dphi}{d _\phi}
\newcommand{\dist}{\mbox{dist}}

\title{Motion by anisotropic mean curvature as sharp interface limit of an inhomogeneous and anisotropic Allen-Cahn equation
}
\author{ }
\date{}

\begin{document}

\maketitle \vspace{-20 mm}

\begin{center}

{\large\bf Matthieu Alfaro }\\[1ex]
D\'epartement de Math\'ematiques CC 051, Universit\'e Montpellier II,\\
Place Eug\`ene Bataillon, 34095 Montpellier Cedex 5, France, \\[2ex]

{\large\bf Harald Garcke}\\[1ex]
Naturwissenschaftliche Fakult\"at I - Mathematik, \\
Universit\"at Regensburg, 93040 Regensburg, Germany, \\[2ex]

{\large\bf Danielle Hilhorst }\\[1ex]
Laboratoire de Math\'ematiques, Analyse Num\'erique et EDP, \\
Universit\'e de Paris Sud, 91405 Orsay Cedex, France, \\[2ex]

{\large\bf Hiroshi Matano}\\[1ex]
Graduate School of Mathematical Sciences, University of Tokyo,\\
3-8-1 Komaba, Tokyo 153-8914, Japan,\\[2ex]

{\large\bf Reiner Sch\"a\-tzle}\\[1ex]
Mathematisches Institut, Arbeitsbereich Analysis,\\
Universit\"at T\"ubingen, Auf der Morgenstelle 10, D-72076
T\"ubingen, Germany.\\[2ex]

\end{center}

\vspace{15pt}

\begin{abstract}
We consider the spatially inhomogeneous and anisotropic
reaction-diffusion equation $u_t=m(x)^{-1}\div[m(x)a_p(x,\n
u)]+\ep ^{-2}f(u)$, involving a small parameter $\ep
>0$ and a bistable nonlinear term whose stable equilibria are $0$
and $1$. We use a Finsler metric related to the anisotropic
diffusion term and work in relative geometry. We prove a weak
comparison principle and perform an analysis of both the
generation and the motion of interfaces. More precisely, we show
that, within the time scale of order $\ep^2|\ln\ep|$, the unique
weak solution $\ue$ develops a steep transition layer that
separates the regions $\{\ue \approx 0\}$ and $\{\ue \approx 1\}$.
Then, on a much slower time scale, the layer starts to propagate.
Consequently, as $\ep \rightarrow 0$, the solution $u^\ep$
converges almost everywhere to $0$ in $\om_t ^-$ and $1$ in $\om
_t ^+$, where $\om_t ^-$ and $\om _t ^+$ are sub-domains of
$\Omega$ separated by an interface $\Gamma _t$, whose motion is
driven by its anisotropic mean curvature. We also prove that the
thickness of the transition layer is of order $\ep$.\\

\end{abstract}

\section{Introduction}\label{intro-aniso}

Evolution laws for interfaces frequently appear in materials
science, differential geometry and image processing. In this paper
we relate so called diffuse and sharp interface models in which
interfaces evolve according to an evolution law, which involves
anisotropic and inhomogeneous driving forces. The evolution
equations we will consider in particular decrease an inhomogeneous
interfacial energy. A diffuse interface model is based on a free
energy which includes gradient terms and, in this paper, the
energy is assumed to be of the following Ginzburg-Landau type
\begin{equation*}
\mathcal{F} (u) = \int_\Omega[a(x,\nabla u) +\frac{1}{\ep^2}
W(u)]m(x)dx\,,
\end{equation*}
where $\Omega\subset\mathbb{R}^N$, $N\ge 2$, is a bounded domain
with smooth boundary, $\ep>0$ is a small parameter related to the
thickness of a diffuse interfacial layer, $W$ is a double well
potential with wells of equal depth and $m$ is a positive
function. We will allow $a$ to be $x$-dependent and anisotropic,
i.e. the value of $a$ will depend on the direction of $\nabla u$.
Taking the gradient flow of $\mathcal{F}$ with respect to the
weighted $L^2$-inner-product $(u,v) = \int_\Omega u(x)v(x)m(x)dx$
leads to the following initial boundary value problem for an
inhomogeneous and anisotropic Allen-Cahn equation
\[(P^\ep) \quad\begin{cases}
u_t=\surm \div \Big[m(x) a_p (x,\n u)\Big]+\edeux  f(u)
&\textrm{ in }\Omega \times (0,+\infty),\vspace{3pt}\\
a_p (x,\n u)\cdot\nu = 0 &\textrm{ on }\partial \Omega \times (0,+\infty),\vspace{3pt}\\
u(x,0)=u_0(x) &\textrm{ in } \Omega,
\end{cases}\]
where $f(u)=-W'(u)$, $\nu$ is the Euclidean unit normal vector
exterior to $\partial\Omega$ and $a_p$ refers to differentiation
with respect to the variable corresponding to $\nabla u$. We
easily derive that solutions to $(P^\ep)$ fulfill
\begin{equation*}
\frac{d}{dt}\mathcal{F}(u) = -\int_\Omega (u_t)^2 m(x)dx \le 0\,,
\end{equation*}
i.e. $\mathcal{F}$ serves as a Lyapunov function.

It can be shown that under appropriate assumptions, see
\cite{bouch}, \cite{owen}, \cite{owenstern}, the energies
$\ep\mathcal{F}$ converge in the sense of $\Gamma$-convergence to
an anisotropic functional defined for hypersurfaces, i.e. in the
limit $\ep\to 0$ the interface is sharp. For a smooth hypersurface
$\Gamma$ the limiting energy becomes
\begin{equation*}
\int_\Gamma\sqrt{2a(x,n(x))}\, m(x)d\mathcal{H}^{N-1}(x)\,,
\end{equation*}
where $n$ is a suitable Euclidean unit normal vector to $\Gamma$
and $d\mathcal{H}^{N-1}$ refers to integration with respect to the
$(N-1)$-dimensional Hausdorff measure. Anisotropic energies for
hypersurfaces can be analyzed in the context of Finsler geometry.
If one considers the steepest descent of the anisotropic surface
energy in relative geometry, where geometric quantities such as
curvature and normal velocity are computed within the context of a
Finsler metric, one obtains an anisotropic and inhomogeneous
generalization of mean curvature flow. In fact the moving
interface $\Gamma_t$ evolves according to the law
\[ (P^0)\quad\begin{cases}
\di{\frac{m(x)}{\sqrt{2a(x,n)}}}V_n=-\div
\Big[\di{\frac{m(x)}{\sqrt{2a(x,n)}}}\; a_p (x,n)\Big]
&\text{ on }\Gamma _t,\vsp\\
\Gamma _t\Big|_{t=0}=\Gamma _0,
\end{cases} \]
where $V_n$ is the normal velocity of $\Gamma _t$. We will show
below that this equation can be rewritten in the relative geometry
associated with a Finsler metric; then it has the form
\[ (P^0)\quad\begin{cases}
V_{n,\phi}=-\kappa _\phi &\text{ on }\Gamma _t, \vsp\\
\Gamma _t\Big|_{t=0}=\Gamma _0,
\end{cases} \]
where $n_\phi$, $V_{n,\phi}$ and $\kappa _\phi$ are, respectively,
the anisotropic unit normal in the exterior direction, the
anisotropic normal velocity of $\Gamma _t$ in the $n_\phi$
direction, and the anisotropic mean curvature at each point of
$\Gamma_t$. In the isotropic homogeneous case one recovers the
mean curvature flow $V_n =-\kappa$. We refer to a paper by
Bellettini and Paolini \cite{BE} and Section \ref{finsler-aniso}
for details.

\vskip 8 pt It is the goal of this paper to rigorously prove that
Problem $(P^\ep)$ converges  to the anisotropic inhomogeneous mean
curvature flow $(P^0)$, as $\ep\to 0$, and to give an optimal
error estimate between the solutions of $(P^\ep)$ and those of
$(P^0)$. We remark that a formal derivation is already contained
in the paper by Bellettini and Paolini \cite{BE}.

Before going into the details, we note that $\Pe$ includes the
following equations as special cases: the spatially inhomogeneous
diffusion equation
\begin{equation}\label{souscas1-aniso}
u_t=\div (A(x) \n u)+\edeux f(u)\,,
\end{equation}
where $A(x)$ is a positive definite symmetric matrix depending on
$x$; the fully anisotropic equation
\begin{equation}\label{souscas2-aniso}
u_t=\div\big(a_p(\n u)\big)+\edeux f(u).
\end{equation}
The significant difference between \eqref{souscas1-aniso} and
\eqref{souscas2-aniso} is that the anisotropy in
\eqref{souscas2-aniso} depends on the solution $u$ itself, while
it does not in \eqref{souscas1-aniso}. In other words, in
\eqref{souscas2-aniso} the dependence of the energy density on the
spatial orientation of the interface can be chosen much more
general when compared with \eqref{souscas1-aniso} where only
ellipsoidal energy densities appear. We refer to Garcke, Nestler
and Stoth \cite{GNS2}, Barrett, Garcke and N\"urnberg
\cite{BGN1,BGN2} for possible anisotropic energy densities. Note
also that we allow $a_{pp}(p)$ to be discontinuous at $p=0$ (see
Remark \ref{physics} below).

We suppose in what follows that $W(u)$ is a double-well potential
with equal well-depth, taking its global minimum value at $u=0$
and $u=1$. More precisely we assume that $f=-W'$ is smooth and has
exactly three zeros $0<a<1$ such that
\begin{equation}\label{der-f-aniso}
f'(0)<0, \quad f'(a)>0, \quad f'(1)<0,
\end{equation}
and that
\begin{equation}\label{int-f-aniso}
\int _ {0} ^ {1} f(u)du=0.
\end{equation}

\begin{rem}\label{unbalanced-remark} Note that we could also
consider the case where $f$ is slightly unbalanced by order $\ep$
so that $\int _0 ^1 f(u)du=O(\ep)$ stands instead of
\eqref{int-f-aniso}. In this case, the singular limit of $\Pe$
will have an additional driving force term in $(P^0)$. See Remark
\ref{second-term-vanishes} for details.\qed
\end{rem}

The assumptions concerning the anisotropic term are the following.
\begin{enumerate}
\item $a(x,p)$ is a real valued function, of class
$C^{3+\vartheta}_{loc}$ (for some $0<\vartheta <1$) on $\ombar
\times \R ^N \backslash \{ 0 \}$;
 \item $a(x,p)$ is positive on $\ombar \times \R ^N \backslash \{ 0 \}$;
 \item $a(x,\cdot)$ is
strictly convex for all $x \in \ombar$;
 \item $a(x,p)$ is homogeneous of degree two in the $p$ variable, i.e.
\begin{equation}\label{a-homogeneite-aniso}
a(x,\alpha p)=\alpha ^2 a(x,p) \quad \text { for all } (x,p) \in
\ombar \times \R ^N \backslash \{ 0 \}, \text { all } \alpha \neq
0.
\end{equation}
\end{enumerate}
If $p$ is given by $p=(p_1,\cdots,p_N)$, the vector valued
function $a_p$ is defined by $a_p(x,p)=\Big(\frac{\partial
a}{\partial p_1}(x,p),\cdots,\frac{\partial a}{\partial
p_N}(x,p)\Big)$, and the matrix valued function $a_{pp}$ by $
a_{pp}(x,p)=\Big(\frac{\partial ^2 a}{\partial p_j
\partial p_i}(x,p)\Big)$. Moreover,
for a vector $p=(p_1,\cdots,p_N)$ and a matrix $A=(a_{ij})$, we
use the notations
$$
|p|=\max _i |p_i| \quad\quad\quad \text{and}\quad\quad
|A|=\max_{i,j}|a_{ij}|.
$$

\vskip 6pt
\begin{rem}\label{der-a-aniso}
The fact that $a$ is homogeneous of degree two implies that, for
all $(x,p) \in \ombar \times \R ^N\backslash \{ 0 \}$, all $\alpha
\neq 0$,
\begin{align*}
a_p(x,\alpha p)&=\alpha a_p(x,p),\;a_{pp}(x,\alpha p)=a_{pp}(x,p),\\
a_p(x,\alpha p)\cdot p&=2\alpha a(x,p),\\
a_{pp}(x,\alpha p) p&=a_p(x,p).
\end{align*}
By setting $a(x,0)=0$ and $a_p(x,0)=0$, one can understand that
$a(x,p)$ is of class $C^1$ on the whole of $\ombar \times \R
^N$.\qed
\end{rem}

\begin{rem}\label{physics}
In many important applications in physics, $a_{pp}(x,p)$ is
discontinuous at $p=0$ and this makes Problem $(P^\ep)$ singularly
parabolic. Because of lack of uniform parabolicity, our analysis
becomes more involved than the case \eqref{souscas1-aniso} or the
case of isotropic Allen-Cahn equation studied in \cite{AHM},
\cite{BS}, \cite{C1}, \cite{MS1, MS2}.\qed
\end{rem}

 We assume that $m:\om \rightarrow (0,+\infty)$ is a
function of class $C^2$ such that $0<m_1\leq m(x) \leq
m_2<+\infty$ for any $x\in\om$, and that $\n m$ and $D^2m$ are in
$L^\infty (\om)$, where $D^2 m (x):=\Big(\di{\frac{\partial ^2
m}{\partial x_j
\partial x_i}(x)}\Big)$.

\begin{rem}\label{megaldet}
Observe that
\begin{equation}\label{mdiv}
\surm \div \Big[m(x) a_p (x,\n u)\Big]=\div a_p(x,\n u)+\n \log
m(x)\cdot a_p(x,\n u)\,,
\end{equation}
so that our equation contains the generalized Allen-Cahn
  equations discussed in Bellettini, Paolini \cite{BE}, Bellettini, Paolini and
Venturini \cite{BE1}.\qed
\end{rem}

We also assume that the initial data $u_0 \in C^2(\ombar)$, and
define $C_0$ as
\begin{equation}\label{int1-aniso}
C_0:=\|u_0\|_{C^0(\ombar)}+\|\n u_0\|_{C^0(\ombar)}+\|
D^2u_0\|_{C^0(\ombar)}.
\end{equation}
Furthermore we define the \lq \lq initial interface" $\Gamma _0$
by
\[\Gamma _0:=\{x\in\om, \; u_0(x)=a \},\]
and suppose that $\Gamma _0$ is a $C^{3+\vartheta}$ closed
hypersurface without boundary ($0<\vartheta<1$), such that, $n$
being the Euclidian unit normal vector exterior to $\Gamma _0$,
\begin{equation}\label{dalltint-aniso}
\Gamma _0 \subset\subset \Omega \quad \mbox { and } \quad \n
u_0(x) \neq 0\quad\text{if $x\in\Gamma _0,$}
\end{equation}
\begin{equation}\label{initial-data-aniso}
u_0>a \quad \text { in } \quad \om ^+ _0,\quad u_0<a \quad \text {
in } \quad \om ^- _0 ,
\end{equation}
where $\om ^- _0$ denotes the region enclosed by $\Gamma _0$ and
$\om ^+ _0$ the region enclosed between $\partial \om$ and $\Gamma
_0$.

For $T>0$, we set $Q_T=\om \times (0,T)$. We define below a notion
of weak solutions of Problem $\Pe$. For this definition, it is
sufficient to only suppose that $u_0 \in H^1(\om)\cap
L^\infty(\om)$.

\begin{defi}\label{definition-weaksol-aniso}
A function $\ue \in L^2(0,T;H^1(\om))\cap L^\infty(Q_T)$ is a weak
solution of Problem $\Pe$, if
\begin{enumerate}
\item $\ue _t \in L^2(Q_T)$,
 \item
$a_p(x,\n \ue (x,t)) \in L^\infty(0,T;L^2(\om))$,
 \item
$\ue(x,0)=u_0(x)$ for almost all $x\in\om$, \item $\ue$ satisfies
the integral equality
\begin{equation}\label{deqexi-aniso}
\int_0^t\int_{\om}\Big[\ue_t \varphi+a_p(x,\n \ue)\cdot \n
\varphi-\edeux f(\ue)\varphi\Big]m(x)dx dt=0\,,
\end{equation}
for all nonnegative function $\varphi\in L^2(0,T;H^1(\om)) \cap
L^\infty(Q_T)$ and for all $t\in(0,T)$.
\end{enumerate}
\end{defi}

\vskip 8 pt One may prove, using monotonicity and compactness
arguments as is done in \cite{B}, \cite{BM}, that Problem $\Pe$
possesses a unique weak solution which we denote by $\ue$. As $\ep
\rightarrow 0$, the qualitative behavior of this solution is the
following. In the very early stage, the anisotropic diffusion term
is negligible compared with the reaction term $\ep ^{-2}f(u)$.
Hence, rescaling time by $\tau=t/\ep^2$, the equation is well
approximated by the ordinary differential equation $u_\tau=f(u)$.
In view of the bistable nature of $f$, $\ue$ quickly approaches
the values $0$ or $1$, the stable equilibria of the ODE, and an
interface is formed between the regions $\{\ue\approx 0\}$ and
$\{\ue\approx 1\}$. Once such an interface has been developed, the
anisotropic diffusion term becomes large near the interface, and
comes to balance with the reaction term so that the interface
starts to propagate, on a much slower time scale.

To understand such interfacial behavior, we have to study the
singular limit of $\Pe$ as $\ep\rightarrow 0$. Then the limit
solution $\tilde u (x,t)$ is a step function taking the values $0$
and $1$ on the sides of the moving interface $\Gamma _t$. In the
case of the usual Allen-Cahn equation, it is well known that
$\Gamma _t$ evolves according to the mean curvature flow
$V_n=-\kappa$ and we will show in this paper that the sharp
interface limit of $(P^\ep)$ is given by $(P^0)$.

Using the theory of analytic semigroups (see e.g. Lunardi
  \cite{Lun}) it is possible to show that the limit
Problem $(P^0)$ possesses locally in time a unique smooth
solution. More precisely, there exists a $T>0$ such that Problem
$(P^0)$ has a unique solution $\Gamma=\bigcup _{0\leq t< T}(
\Gamma _t\times\{ t\})$ which satisfies $\Gamma \in
C^{3+\vartheta,(3+\vartheta)/2}$. For proofs of the local in time
existence of solutions of related limit problems, we also refer
the reader to \cite{GG} and the discussion at the end of
  Chapter 1 in \cite{giga}.

For each $t\in (0,T)$, we define $\om^-_t$ as the region enclosed
by the hypersurface $\Gamma _t$ and $\om^+_t$ as the region lying
between $\partial \om$ and $\Gamma _t$. Furthermore we define a
step function $\tilde u(x,t)$ by
\begin{equation}\label{u-aniso}
\tilde u(x,t)=\begin{cases} 1 &\text{in } \om^+_t\\
0 &\text{in } \om^-_t \end{cases} \quad\text{for } t\in(0,T).
\end{equation}

It is convenient to present our main result, Theorem
\ref{width-aniso}, in the form of a convergence theorem, mixing
generation and propagation. It describes the profile of the
solution after a very short initial period. It asserts that: given
a virtually arbitrary initial data $u_0$, the solution $\ue$
quickly becomes close to $0$ or $1$, except in a small
neighborhood of the initial interface $\Gamma _0$, creating a
steep transition layer around $\Gamma _0$ ({\it generation of
interface}). The time $t^\ep$ for the generation of interface is
of order $\ep ^2|\ln \ep|$. The theorem then states that the
solution $\ue$ remains close to the step function $\tilde u$ on
the time interval $(t^\ep,T)$ ({\it motion of interface}).
Moreover, as is clear from the estimates in the theorem, the
thickness of the transition layer is of order $\ep$.

\begin{thm}[Generation, motion and thickness of interface]\label{width-aniso}
Let $\eta$ be an arbitrary constant satisfying $0< \eta <
\min(a,1-a)$ and set
$$
\mu=f'(a).
$$
Then there exist positive constants $\ep _0 $ and $C$ such that,
for all $\,\ep \in (0,\ep _0)$ and for almost all $(x,t)$ such
that $\,t^\ep \leq t \leq T$, where
\begin{equation}\label{tep-aniso}
t^\ep:=\mu ^{-1}  \ep ^2 |\ln \ep|\,,
\end{equation}
we have,
\begin{equation}\label{resultat-aniso}
\ue(x,t) \in
\begin{cases}
\,[-\eta,1+\eta]&\quad\text{if}\quad
x\in\mathcal N_{C\ep}(\Gamma_t)\vsp\\
\,[-\eta,\eta]&\quad\text{if}\quad x\in\om_t^-\setminus\mathcal
N_{C\ep}(\Gamma
_t)\vsp\\
\,[1-\eta,1+\eta]&\quad\text{if}\quad x\in\om
_t^+\setminus\mathcal N_{C\ep}(\Gamma _t),
\end{cases}
\end{equation}
where $\mathcal N _r(\Gamma _t):=\{x\in \om, dist _\phi(x,\Gamma
_t)<r\}$ denotes the $r$-neighborhood of $\Gamma _t$; by $dist
_\phi(x,\Gamma_t)$, we mean the $\delta _\phi$ distance to the set
$\Gamma _t$, where $\delta _\phi$ is the distance associated to a
Finsler metric, whose definition is to be given in Section
\ref{finsler-aniso}.
\end{thm}

\begin{cor}[Convergence]\label{total-aniso}
As $\ep\to 0$, the solution $\ue$ converges to $\tilde u$ almost
everywhere in $\bigcup _{0<t< T}(\om^\pm_t\times\{ t\})$.
\end{cor}

\vskip 8 pt The organization of this paper is as follows. In
Section \ref{finsler-aniso}, we recall notations and results
concerning Finsler metrics that give a natural and efficient
framework for dealing with anisotropic problems. In Section
\ref{formal-aniso}, we perform formal asymptotic expansions in
order to derive the equation for the motion of interface, and
collect useful estimates on stationary solutions of related
problems. In Section \ref{comparison-aniso} we prove a weak
comparison principle for Problem $\Pe$. Such a comparison
principle is rather standard, but, in view of the fact that
$a_{pp}(x,p)$ does not exist at $p=0$, we give a complete proof
for the convenience of the reader. In Section
\ref{generation-aniso}, we prove results on the generation of
interface. For the study of this early time range we construct
sub- and super-solutions by modifying the solution of the ordinary
differential equation $u_t=\ep^{-2}f(u)$. In Section
\ref{motion-aniso}, we construct another pair of sub- and
super-solutions by using the first two terms of the formal
asymptotic expansion given in Section \ref{formal-aniso}. They are
used to study the motion of interface in the later stage. In
Section \ref{s:proof-aniso}, by fitting these two pairs of sub-
and super-solutions together, we prove our main results for $\Pe$:
Theorem \ref{width-aniso} and its corollary.

Let us mention  some earlier works on anisotropic problems related
to $(P^\ep)$. In \cite{BFP}, Bellettini, Colli Franzone and
Paolini study a problem that is slightly more general than $\Pe$
--- by allowing $f$ to be unbalanced in the same way as in Remarks
\ref{unbalanced-remark}, \ref{second-term-vanishes} of the present
paper --- and derive a very fine error estimate between the formal
asymptotic and actual solutions of $\Pe$. We also refer to the
articles \cite{ESch, ESch2}, by Elliott and Sch\"a\-tzle, on a
similar but slightly different problem where the potential $W(u)$
is a double obstacle type, namely $W(u)=+\infty$ for $u\notin
(0,1)$. For the spatially homogeneous case $a(x,p)=a(p)$ they
prove convergence of the anisotropic diffusion problem to an
anisotropic curvature flow similar to $(P^0)$. Note that their
second paper \cite{ESch2} considers a kinetic term of the form
$\beta(\n u)u_t$, which makes the meaning of solutions very weak,
hence they are treated in the framework of viscosity solutions.

However, in these papers, the authors consider only a very
restricted class of initial data, namely those having a specific
profile with well-developed transition layer. More precisely they
prove that if the initial data is very close to the typical
profile that appears in the formal asymptotic expansions of the
moving interface, then the solution remains close to the formal
asymptotic for $0\leq t \leq T$. In other words the generation of
interface from arbitrary initial data is not studied there.
Summarizing, they have obtained a very fine error estimate --- of
order $O(\ep^2)$ or higher --- between the solutions of specific
initial data and formal asymptotic, while, in the present paper,
we consider convergence of solutions of $\Pe$ with virtually
arbitrary initial data to solutions of $(P^0)$, with an error
estimate of order $O(\ep)$. Therefore, the two results are both
for the convergence of $\Pe$ to $(P^0)$, but they are of different
nature. Note that, as far as the thickness of the interface is
concerned, our $O(\ep)$ estimate is optimal (see \cite{AHM} for
details).

In \cite{BHW}, Bene{\v{s}}, Hilhorst and Weidenfeld study both the
generation and the motion of interface for an anisotropic
Allen-Cahn equation which is related to ours. Nevertheless, their
equation is slightly less general since they do not allow
$x$-dependence in $a(x,p)$. Moreover, with their sub- and
super-solutions, they cannot achieve the optimal $O(\ep)$ estimate
of the thickness of the interface.

For numerical simulations for problems $(P^\ep)$ and $(P^0)$ we
refer to Bene{\v{s}}, Mikula \cite{BM}, Garcke, Nestler, Stoth
\cite{GNS2}, Barrett, Garcke, N\"urnberg \cite{BGN1, BGN2} and
Paolini \cite{P}.

\section{Finsler metrics and the anisotropic context}\label{finsler-aniso}

In this section we explain the technique of Bellettini, Paolini
\cite{BE}, Bellettini, Paolini and Venturini \cite{BE1} to apply
Finsler metric to analyze anisotropic nonlinear problems. The idea
is to endow $\mathbb{R}^N$ with the distance obtained by
integrating the Finsler metric which makes otherwise lengthy
computations remarkably simpler. For the convenience of the
reader, we first recall basic properties of Finsler metrics as
stated in \cite{BE}, \cite{BE1}. For more details and proofs, see
these references.

\subsection{Finsler metrics}

Suppose that $\phi:\om \times \mathbb{R}^N \rightarrow
[0,+\infty)$ is a continuous function satisfying the properties
\begin{align}
\phi (x,\alpha\xi)=|\alpha|\phi (x,\xi) &\quad\ \textrm{for all}\
\; (x,\xi) \in \om \times \mathbb{R}^N \quad \text{and all} \; \
\alpha \in \R,\vspace{3pt}\label{fins1-aniso}\\
\lambda _0|\xi| \leq \phi (x,\xi) \leq \Lambda _0|\xi| &\quad\
\textrm{for all}\ \; (x,\xi) \in \om \times
\mathbb{R}^N,\label{fins2-aniso}
\end{align}
for two suitable constants $0<\lambda _0\leq \Lambda _0<+\infty$.
We say that $\phi$ is strictly convex if, for any $x \in \om$, the
map $\xi\mapsto \phi ^2 (x,\xi)$ is strictly convex on
$\mathbb{R}^N$. We shall indicate by
$$
B_\phi(x)=\{\xi \in \mathbb{R}^N, \phi(x,\xi)\leq 1\}
$$
the unit sphere of $\phi$ at $x \in \om$.

The dual function $\phi ^0:\om \times \mathbb{R}^N \rightarrow
[0,+\infty)$ of $\phi$ is defined by
\begin{equation}\label{dual-aniso}
\phi ^0 (x,\xi^*)=\sup \big{\{}\xi^* \cdot \xi,\; \xi\in
B_\phi(x)\big{\}}\,,
\end{equation}
for any $(x,\xi) \in \om \times \mathbb{R}^N$. One can prove that
$\phi ^0 $ is continuous, convex, satisfies properties
\eqref{fins1-aniso} and \eqref{fins2-aniso}, and that $\phi
^{00}$, the dual function of $\phi ^0$, coincides with the convex
envelope of $\phi$ with respect to $\xi$.

We say that $\phi$ is a (strictly convex smooth) Finsler metric,
and we shall write $\phi \in \mathcal M(\om)$ if, in addition to
properties \eqref{fins1-aniso} and \eqref{fins2-aniso}, $\phi$ and
$\phi ^0$ are strictly convex and of class $C ^2$ on $\om \times
\mathbb R ^N \setminus \{0\}$. In particular $\phi ^{00}=\phi$.

We denote by $\delta _\phi$ the integrated distance associated to
$\phi \in \mathcal M (\om)$, that is, for any $(x,y) \in \om$, we
set
\begin{equation}\label{dist-phi-aniso}
\delta _\phi (x,y)=\inf \Big{\{} \int _0 ^1 \phi(\gamma
(t),\dot{\gamma }(t))dt\,;\ \gamma \in W^{1,1}\big([0,1];\om\big),
\gamma(0)=x, \gamma (1)=y \Big{\}}.
\end{equation}
In the special case of the Euclidian metric, the function $\phi$
is given by $\phi(x,p)=\phi(p)=({p_1} ^2+\cdots +{p_N} ^2)^{1/2}$,
so that $\delta _\phi$ reduces to the usual distance.

Given $\phi \in \mathcal M (\om)$ and $x \in \om$, let
$T^0(x,\cdot):\mathbb R^N \to \mathbb R^N$ be the map defined by
\begin{equation}\label{T0-aniso}T^0(x,\xi^*)=\begin{cases} \phi ^0(x,\xi ^*)
\phi ^0 _ p (x,\xi^*) &\text{if } \xi ^* \in \mathbb R ^N \setminus \{0\}\\
0 &\text{if } \xi ^* =0.\end{cases}
\end{equation}
Here $\phi ^0 _p$ denotes the gradient with respect to $p$
whenever we regard $\phi ^0(x,p)$ as a function of two variables
$x$ and $p$.

If $u:\om \to \mathbb R$ is a smooth function with non-vanishing
gradient, we define the anisotropic gradient by
\begin{equation}\label{nabla-phi-aniso}
\n _\phi u =T^0(x,\n u )=\phi ^0(x,\n u) \phi ^0 _ p (x,\n u).
\end{equation}
If $\eta:\om \to \mathbb R ^N$ is a smooth vector field, we define
the $m$-divergence operator by
\begin{equation}\label{div-phi-aniso}
\div _ m \eta =\surm \div \big[m(x)\eta\big]=\div \eta+ \n \log
m(x)\cdot \eta\,,
\end{equation}
and then the $m$-anisotropic Laplacian by
\begin{equation}
\Delta _{\phi,m} u =\div _m \n _\phi u.
\end{equation}
Note that in \cite{BE}, \cite{BE1} $m$ is related to $\phi$ while
in the present paper $m$ is a given function independent of
$\phi$. Nonetheless, in the sequel, we shall use the simpler
notation $\Delta _\phi :=\Delta _{\phi,m}$.

As in the isotropic case, if $\Gamma _t$ is a smooth hypersurface
of $\om$ at time $t$, and $n$ the outer normal vector to $\Gamma
_t$ (in the Euclidian sense), we define $n_ \phi$ the
$\phi$-normal vector to $\Gamma _t$ and $\kappa _\phi$ the
$\phi$-mean curvature of $\Gamma _t$ by
\begin{equation}
n_{\phi}=\phi ^0_p(x,n),\quad \kappa_{\phi}=\div _m n_{\phi}.
\end{equation}

Furthermore, if $\psi$ is a smooth function with non-vanishing
gradient such that $\Gamma _t=\{x\in\om,\; \psi(x,t)=0\}$, and
$\psi$ is positive outside $\Gamma _t$ and negative inside, then
\begin{align}\label{phi-normal-aniso}
n&=\di{\frac{\n\psi}{|\n\psi|}},\quad\quad\quad\quad\quad\quad n_{\phi}=\phi ^0_p(x,\n\psi),\\
\label{phi-curvature-aniso} \kappa&= \div
\di{\frac{\n\psi}{|\n\psi|}},\quad\quad\quad\quad\kappa_{\phi}=\div
_m \phi ^0 _p(x,\n\psi),
\end{align}
on $\Gamma _t$. We also define the normal velocity of $\Gamma _t$
and the $\phi$-normal velocity of $\Gamma _t$ by
\begin{equation}\label{phi-velocity-aniso}
V_n =- \di{\frac{\psi_t}{|\n\psi|}},\quad\quad\quad\quad
V_{n,\phi} = -\di{\frac{\psi_t}{\phi ^0(x,\n\psi)}}.
\end{equation}

To conclude these preliminaries, we quote a theorem proved in
\cite{BE1}.
\begin{thm}\label{egal1-aniso}
Let $\om$ be connected, and let $\phi \in \mathcal M (\om)$. Let
$\delta _\phi$ be the integrated distance associated to $\phi$.
Let $C\subseteq \om$ be a closed set, and let $\dist  _\phi(x,C)$
be the $\delta _ \phi$ distance to the set $C$ defined by
\begin{equation}\label{dist-belletini-aniso}
\mbox{dist}_\phi(x,C)=\inf \big{\{}\delta _ \phi (x,y)\;,\; y\in
C\big{\}}.
\end{equation}
Then
\begin{equation}\label{bellettini-general-aniso}
\phi ^0 \big(x, \n \mbox{dist}_\phi(x,C)\big)=1\,,
\end{equation}
at each point $x \in \om \setminus C$ where
$\mbox{dist}_\phi(\cdot,C)$ is differentiable.
\end{thm}

\vskip 8 pt In the special case of the Euclidian metric, note that
\eqref{bellettini-general-aniso} reduces to the property that $|\n
d|=1$.

\subsection{Application to the anisotropic Allen-Cahn
equation}\label{sub-appl-AAC}

We set, for all $(x,p)\in\om\times\R ^N$,
\begin{equation}\label{case-aniso}
\phi ^0(x,p)=\sqrt{2a(x,p)}.
\end{equation}
First, since $a(x,\cdot)$ is homogeneous of degree two, $\phi ^0$
satisfies assumptions \eqref{fins1-aniso} and \eqref{fins2-aniso}
with the constants
\begin{equation}\label{lambda0-aniso}
\lambda_0=[2\min _{x \in \ombar, |p|=1} a(x,p)]^{1/2}>0 \quad
\text {and}\quad \Lambda _0=[2\max _{x \in \ombar, |p|=1}
a(x,p)]^{1/2}>0.
\end{equation}
By the hypotheses on $a(x,p)$, we see that $\phi ^0$ is strictly
convex and of class $C^2$ on $\om \times \mathbb R ^N \setminus
\{0\}$; moreover, by Remark \ref{der-a-aniso}, $\phi ^0$ is
continuous on the whole of $\om \times \mathbb{R}^N $. It follows
that $\phi$ is a Finsler metric and the above theory applies. We
have
\begin{equation}\label{T0nous-aniso}
T^0(x,p)=\begin{cases}a_p(x,p) &\text{if }
p \in \mathbb R ^N \setminus \{0\}\\
\;\;\;\;\;\;\; 0 &\text{if } p=0.\end{cases}
\end{equation}
Let $\Gamma=\bigcup _{0\leq t < T}(\Gamma _t\times\{t\})$ be the
unique solution of the limit geometric motion Problem $(P^0)$ and
let $\widetilde d$ be the signed distance function to $\Gamma$
defined by
\begin{equation}\label{eq:distance-aniso}
\widetilde d (x,t)=
\begin{cases}
\phantom{-}\mbox{dist}(x,\Gamma _t) \quad \text{for }x\in\om _t^+, \\
-\mbox{dist}(x,\Gamma _t) \quad \text{for } x\in\om _t^- ,
\end{cases}
\end{equation}
where $\mbox{dist}(x,\Gamma _t)$ is the distance from $x$ to the
hypersurface  $\Gamma _t$ in $\om$. Let $\widetilde d_\phi$ be the
anisotropic signed distance function to $\Gamma$ defined by
\begin{equation}\label{eq:dist-aniso}
\widetilde d_\phi (x,t)=
\begin{cases}
\phantom{-}\mbox{dist}_\phi(x,\Gamma _t) \quad \text{for }x\in\om _t^+,\\
-\mbox{dist}_\phi(x,\Gamma _t) \quad \text{for } x\in\om _t^- ,
\end{cases}
\end{equation}
where $\mbox{dist}_\phi(x,\Gamma _t)$ denotes the $\delta _ \phi$
distance to the set $\Gamma _t$ defined in
\eqref{dist-belletini-aniso}. By Theorem \ref{egal1-aniso}, the
following equality holds
\begin{equation}\label{bellettini-aniso} 2a(x,\n \widetilde
d_\phi(x,t))=1 \quad \quad \text{ in a neighborhood of } \Gamma
_t.
\end{equation}
By setting $\psi=\widetilde d$ and $\psi=\widetilde d _\phi$ in
the second equalities in \eqref{phi-normal-aniso},
\eqref{phi-curvature-aniso}, \eqref{phi-velocity-aniso}, we obtain
two equivalent expressions of the $\phi$-normal vector, the
$\phi$-mean curvature and the $\phi$-normal velocity:
\begin{align}
\label{link1-aniso}&n_ \phi=\frac{1}{\sqrt{2a(x,\n \widetilde
d)}}\; a_p(x,\n \widetilde d) = \; a_p(x,\n
\widetilde d _\phi),\vsp\\
\label{link2-aniso}&\kappa _ \phi=\div
_m\Bigg[\frac{1}{\sqrt{2a(x,\n \widetilde d)}}\;a_p(x,\n
\widetilde d)\Bigg]=\div _m\Big[a_p(x,\n
\widetilde d _\phi)\Big],\vsp\\
\label{link3-aniso}&V_{n,\phi}=-\,\frac{1}{\sqrt{2a(x,\n
\widetilde d)}}\; \widetilde d _t =-(\widetilde d _\phi)_t.
\end{align}

\vskip 8 pt The end of this section is devoted to the anisotropic
Laplacian
\begin{align}
\Delta _{\phi} u &=\surm \div\big[m(x)a_p(x,\n u)\big]\vsp\label{laplace-aniso1}\\
&=\div a_p(x,\n u)+\n \log m(x)\cdot a_p(x,\n
u)\label{laplace-aniso2}.
\end{align}
In the case of Finsler metrics, it turns out that the term $\Delta
_{\phi} u$ may be less regular than $\Delta u$. Nevertheless, we
show below a boundedness property of the anisotropic Laplacian
(see \cite{BHW} for a related property).

\begin{lem}\label{aniso-lapla-aniso}
There exists a positive constant $C_L$ such that, for all $u \in
C^{2,1}(\Q)$, the following inequality holds:
\begin{equation}\label{aniso-lapla-eq-aniso}
|\Delta _\phi u(x,t)|\leq C_L(| \n u (x,t)|+| D^2u (x,t)|)
\quad\quad \text{for all $(x,t) \in \Q$.}
\end{equation}
\end{lem}

{\noindent \bf Proof.} In view of \eqref{laplace-aniso2}, it is
sufficient to deal with the term $\div a_p(x,\n u)$. We can, with
no loss of generality, ignore the dependence on time.

First, assume that $x$ is such that $\n u(x)\neq 0$. Regarding
$a(x,p)$ as a function of two variables $x$ and
$p=(p_1,\cdots,p_n)$, we obtain, by a straightforward calculation,
that
\begin{equation}\label{rennes-aniso}
\div  a_p(x,\n u(x))=\sum _j \frac{\partial ^2 a}{\partial x_j
\partial p_j}(x,\n u(x))+\sum _{i,j}\frac{\partial ^2 a}{\partial p_i
\partial p_j}(x,\n u(x))\frac{\partial ^2 u}{\partial x_i
\partial x_j}(x).
\end{equation}
It follows from the homogeneity properties that
\begin{align*}
|\div a_p(x,\n u(x))| \leq &|\n u(x)| \sum _j \max _{y\in\ombar,
|p|=1} \Big|\frac{\partial ^2 a}{\partial x_j
\partial p_j}(y,p) \Big|\\
&+|D^2 u (x)| \sum _{i,j}\max _{y \in \ombar,
|p|=1}\Big|\frac{\partial ^2 a}{\partial p_i
\partial p_j}(y,p)\Big|\,,
\end{align*}
where we have used that $a$ is of class $C^2$ on the compact set
$\ombar \times \{|p|=1\}$. This proves
\eqref{aniso-lapla-eq-aniso} under the assumption $\n u(x) \neq
0$.

Now assume that $x$ is such that $\n u(x)=0$. We have to proceed
in a slightly different way since $a_{pp}(x,0)$ does not make
sense. The operator $a_p(x,\cdot)$ is homogeneous of degree one so
that, for any direction $\zeta$,
$$
t^{-1}(a_p(x,t\zeta)-a_p(x,0))=a_p(x,\zeta).
$$
We denote by $(e_1,\cdots,e_N)$ the Euclidian basis of
$\mathbb{R}^N$. It follows from the above equality that
$a_p(x,\cdot)$ admits at the point $0$ partial derivatives in any
direction $e_i$ and
\begin{equation}\label{direction-derivative-aniso}
\frac{\partial a_p(x,\cdot)}{\partial p_i} (0)=a_p(x,e_i)\,,
\end{equation}
which, in turn, implies that
\begin{equation}\label{avantveille-aniso}
\frac{\partial}{\partial p_i}\frac{\partial a}{\partial
p_j}(x,0)=\frac{\partial a}{\partial p_j}(x,e_i).
\end{equation}
Note that, since $a_p(x,\cdot)$ is homogeneous of degree one, the
first term in \eqref{rennes-aniso} vanishes at the point $(x,0)$.
It follows from \eqref{rennes-aniso} and \eqref{avantveille-aniso}
that, in the case where $\n u(x)=0$,
\begin{align*}
|\div a_p(x,\n u(x))|&=\Big|\sum _{i,j}\frac{\partial a}{\partial
p_j}(x,e_i)\frac{\partial ^2 u}{\partial x_i \partial
x_j}(x)\Big|\\
&\leq | D^2 u(x)| \sum _{i,j}\max _{y \in \ombar, |p|=1}
\Big|\frac{\partial a}{\partial p_j}(y,p)\Big|\,,
\end{align*}
which proves \eqref{aniso-lapla-eq-aniso} in this case as well.
\qed

\section{Formal derivation of the interface motion
equation}\label{formal-aniso}

In this section we derive the equation of interface motion
corresponding to Problem $\Pe$ by using a formal asymptotic
expansion. The resulting interface equation can be regarded as the
singular limit of $\Pe$ as $\ep \to 0$. Our argument goes
basically along the same lines with the formal derivation given by
Nakamura, Matano, Hilhorst and Sch\"atzle \cite{NMHS}: the first
two terms of the asymptotic expansion determine the interface
equation. Though our analysis in this section is for the most part
formal, the results we obtain will help the rigorous analysis in
later sections.\\

Let $u^\ep$ be the solution of $\Pe$. Let $\Gamma=\bigcup _{0\leq
t < T} \Gamma_t\times\{t\} $ be the solution of the limit
geometric motion problem and $\widetilde d_\phi$ the anisotropic
signed distance function to $\Gamma$ defined in
\eqref{eq:dist-aniso}. We then define
\begin{equation*}
Q^+_T = \bigcup_{\,0<t< T}\Omega^+_t \times\{t\},\qquad Q^-_T =
\bigcup_{\,0<t<T}\Omega^-_t \times\{t\}.
\end{equation*}

We also assume that the solution $\ue$ has --- one the one hand
--- the outer expansions (away from the interface $\Gamma$),
\begin{equation} \label{outer-aniso}
\ue(x,t)= \tilde u (x,t) + \ep u_1^{\pm} (x,t) + \ep^2
u_2^{\pm}(x,t) + \cdots \quad\ \hbox{in}\ \ Q^{\pm}_T,
\end{equation}
where $\tilde u$ is the step function defined in \eqref{u-aniso},
and --- on the other hand --- the inner expansion (near $\Gamma$)
\begin{equation}\label{inner-aniso}
\ue(x,t)=\U(x,t,\xi)+\ep U_1 (x,t,\xi)+\ep^2 U_2 (x,t,\xi)+\cdots,
\end{equation}
near $\Gamma$ (the inner expansion), where $U_j(x,t,z)$,
$j=0,1,2,\cdots$, are defined for $x\in \overline \Omega$, $t\geq
0$, $z\in \mathbb{R}$. The stretched space variable
$\xi:=\widetilde d_\phi(x,t)/\ep$ gives exactly the right spatial
scaling to describe the sharp transition between the regions
$\{\ue \approx 0\}$ and $\{\ue \approx 1\}$. We normalize $U_0$ in
such a way that
$$
\U(x,t,0)=a
$$
(normalization conditions). To make the inner and outer expansions
consistent, we require that
\begin{equation}\label{match}
\begin{array}{ll}
\U(x,t,+\infty)=1 , \quad &U_k(x,t,+\infty)= u_k^{+}(x,t),
\vspace{3pt}\\
\U(x,t,-\infty)=0 , \quad &U_k(x,t,-\infty)= u_k^{-}(x,t),
\end{array}
\end{equation}
for all $k \geq 1$ (matching conditions).

In what follows we will substitute the inner expansion
\eqref{inner-aniso} into the parabolic equation in Problem
$(P^\ep)$ and collect the $\ep^{-2}$ and $\ep^{-1}$ terms. For
this purpose, note that if $V=V(x,t,z)$ and $v(x,t)=V(x,t,\xi)$
are real valued functions then $\n v=\ep ^{-1}V_z\n \dtip +\n _x
V$ and $v_t=\ep^{-1}(\dtip) _t V_z +V_t$; if $v$ and $V$ are
vector valued functions we obtain $\div v=\ep^{-1} \n \dtip \cdot
V_z+ \div _x  V$. In the following, we shall use the properties
stated in Remark \ref{der-a-aniso}. A straightforward computation
yields
\begin{align*}\
\ue_t&= \frac 1 \ep (\dtip) _t \UU+ U_{0t}+(\dtip) _t \VV+\ep
U_{1t}+\cdots \\
\nabla  \ue&= \frac 1 \ep \UU \n \dtip +\n _x\U+\VV \n \dtip +\ep
\n _x U_1 +\cdots  \\
a_p(x,\n \ue)&=\frac 1 \ep a_p(x, \UU \n \dtip+\ep \n _x\U+\ep \VV
\n \dtip
+\ep ^2\n _x U_1+\cdots)\\
&=\frac 1 \ep a_p(x, \UU \n \dtip)+a_{pp}(x,\UU \n \dtip)
(\n _x\U+\VV \n \dtip  )+\cdots\\
&=\frac 1 \ep \UU a_p(x,\n \dtip)+a_{pp}(x,\n \dtip) (\n _x\U+\VV
\n \dtip)+\cdots.
\end{align*}
It follows that
\begin{align*}
\n \log m(x)\cdot a_p(x,\n \ue)=\frac 1 \ep \UU \n \log m(x)\cdot
a_p(x,\n \dtip)+\cdots
\end{align*}
and that
\begin{align*}
\div a_p(x,\n \ue)=&\frac 1 \ep \n \dtip \cdot \partial _ z
(a_p(x,\n\ue))+ \div _x (a_p(x,\n \ue))\\
=& \frac {\n \dtip} {\ep}\cdot \Big[\frac {\UUU}{\ep} a_p(x,\n
\dtip)+a_{pp}(x,\n \dtip)(\n _x \UU+\VVV \n \dtip )\Big]\\
&\quad\quad+\frac 1 \ep \big[\n _x \UU \cdot a_p(x,\n\dtip)+\UU
\div a_p(x,\n\dtip)\big]+\cdots\\
=&\frac {1}{\ep ^2} \UUU 2 a(x,\n \dtip)+\frac 1 \ep\Big[2 a(x,\n
\dtip)\VVV \\&\quad\quad +2 \n _x \UU\cdot a_p(x,\n\dtip)+\UU \div
a_p(x,\n \dtip)\Big]+\cdots,
\end{align*}
where the functions $U_i$ $(i = 0, 1)$, as well as their
derivatives, are taken at the point $(x,t,\dtip (x,t)/\ep)$.
Hence, in view of \eqref{bellettini-aniso}, we obtain
\begin{align*}
\div a_p(x,\n \ue)=& \frac {1}{\ep ^2} \UUU+\frac 1 \ep[\VVV
+2 \n _x \UU \cdot a_p(x,\n \dtip)\\
&\quad\quad\quad\quad\quad\quad\quad
 +\UU \div a_p(x,\n \dtip)]+\cdots.
\end{align*}
We also use the expansion
$$
f(\ue)=f(U_0)+\ep U_1f'(U_0)+\cdots.
$$
Next, we substitute the above expressions in the partial
differential equation in Problem $\Pe$. Collecting the $\ep ^{-2}$
terms yields
\begin{equation}\label{eq-phi-aniso}
\UUU+f(\U)=0.
\end{equation}
In view of the normalization and matching conditions, we can now
assert that $\U(x,t,z)=U_0(z)$, where $U_0$ is the unique solution
of the one-dimensional stationary problem
\begin{equation}\label{eq-phi-aniso-pb}
\left\{\begin{array}{ll}
{U_0} '' +f(U_0)=0,\vsp\\
U_0(-\infty)= 0,\quad U_0(0)=a,\quad U_0(+\infty)=1.
\end{array} \right.
\end{equation}
This solution represents the first approximation of the profile of
a transition layer around the interface observed in the stretched
coordinates. We recall standard estimates on $U_0$.
\begin{lem}\label{est-phi-aniso}
There exist positive constants $C$ and $\lambda$ such that
\begin{align*}
0 <1-U_0(z)&\leq Ce^{-\lambda|z|} \quad \text{ for } z\geq 0,\\
0 <U_0(z)&\leq Ce^{-\lambda|z|} \quad \text{ for } z\leq 0.
\end{align*}
In addition to this ${U_0}'>0$ and, for all $j=1,2$,
$$
|D^j U_0(z)|\leq Ce^{-\lambda|z|} \quad \text{ for } z\in
\mathbb{R}.
$$
\end{lem}

\vskip 8 pt Since $U_0$ depends only on the variable $z$, we have
$\n _x {U_0}'=0$. Then, by collecting the $\ep ^{-1}$ terms, we
obtain
\begin{equation}\label{eq-U1-interm-aniso}
\VVV + f'(\U) U_1 =(\dtip) _t {U_0}'- \Delta _\phi \dtip
\;{U_0}'\,,
\end{equation}
which can be seen as a linearized problem for
\eqref{eq-phi-aniso}. The solvability condition for the above
equation, which is a variant of the Fredholm alternative, plays
the key role in deriving the equation of interface motion. It is
given by
\begin{align*}
\int_\R \Big [(\dtip) _t(x,t)-\Delta _\phi \dtip
(x,t)\Big]{{U_0}'}^2(z)dz=0\,,
\end{align*}
for all $(x,t) \in Q_T$. It follows that $(\dtip) _t=\Delta _\phi
\dtip$. In virtue of subsection \ref{sub-appl-AAC}, this equation,
written in relative geometry, reads as
\begin{equation}\label{interface-relativegeo-aniso}
V_{n,\phi}=-\kappa _\phi \quad\ \text{on}\ \; \Gamma _t,
\end{equation}
that is the interface motion equation $(P^0)$, whereas, in the
Euclidian geometry, the same equation reads as
\begin{equation}\label{interface-aniso}
\di{\frac{m(x)}{\sqrt{2a(x,n)}}}\; V_n=-\div
\Big[\di{\frac{m(x)}{\sqrt{2a(x,n)}}}\; a_p (x,n)\Big]\quad\
\text{on}\ \; \Gamma _t.
\end{equation}
Summarizing, under the assumption that the solution $\ue$ of
Problem $\Pe$ satisfies
\begin{equation*}
\ue\to
\begin{cases}
1 &\quad  \textrm{ in } Q_T^+ \\
0 &\quad  \textrm{ in } Q_T^-
\end{cases}\qquad\hbox{as}\ \ \ep\to 0,
\end{equation*}
we have formally proved that the boundary $\Gamma _t$ between
$\Omega_t^-$ and $\Omega_t^+$ moves according to the law
\eqref{interface-relativegeo-aniso} or \eqref{interface-aniso}.

\begin{rem}\label{second-term-vanishes}
To conclude this section, note that \eqref{eq-U1-interm-aniso} now
yields $U_1=0$. In fact, the second term of the asymptotic
expansion vanishes because the two stable zeros of the
nonlinearity $f$ have \lq \lq balanced" stability, or more
precisely because of the assumption $\int _0 ^1 f(u)du=0$. If we
perturb the nonlinearity by order $\ep$, say $f(u)\longleftarrow
f(u)-\ep g(x,t,u)$, the equation in the free boundary problem
contains an additional driving force term and $U_1$ no longer
vanishes. More precisely, the equation will read as
\[
 \, V_{n,\phi}=-\kappa _\phi+ c_0 \int _0 ^1 g(x,t,r)dr
 \quad \text { on } \Gamma_t,
\]
with $c_0$ a constant explicitly determined by the nonlinearity
$f$. We refer to \cite{AHM} for details.\qed
\end{rem}

\section{A comparison principle}\label{comparison-aniso}

This section is devoted to a comparison principle for weak
solutions of Problem $\Pe$. Such a result is rather standard (see
\cite{BHW}), but, since the problem is non-regular where $\n u=0$,
we prove it here
for the self-containedness of the paper.\\

To begin with, we define a notion of sub- and super-solution of
Problem $\Pe$.

\begin{defi}\label{definition-aniso}
A function $\upl \in L^2(0,T;H^1(\om))\cap L^\infty(Q_T)$ is a
weak super-solution of Problem $\Pe$, if
\begin{enumerate}
\item $(\upl)_t \in L^2(Q_T)$,
 \item $\nabla_\phi \upl (x,t)=a_p(x,\n
\upl(x,t)) \in L^\infty(0,T;L^2(\om))$,
 \item $\ue$ satisfies the integral inequality
\begin{equation}\label{deqexi-super-aniso}
\int_0^t\int_{\om}\Big[(\upl)_t \varphi+a_p(x,\n \upl)\cdot \n
\varphi-\edeux f(\upl)\varphi\Big]m(x)dxdt\geq 0,
\end{equation}
for all nonnegative function $\varphi\in L^2(0,T;H^1(\om)) \cap
L^\infty(Q_T)$ and for all $t\in(0,T)$.
\end{enumerate}

We define a weak sub-solution $\um$ in a similar way, by changing
$\geq$ in (\ref{deqexi-super-aniso}) by $\leq$.
\end{defi}

\vskip 8 pt The following remark will be useful when constructing
smooth sub- and super-solutions in later sections.

\begin{rem}\label{rem-sup-aniso}
If $\upl \in C^{2,1}(\Q)$, it is not difficult to see that $\upl$
is a super-solution in the sense defined above if and only if
\begin{enumerate}
\item $a_p(x,\n \upl)\cdot\nu\geq 0$\quad on $\partial\Omega
\times(0,T)$, \item ${\cal L}_0\upl\geq 0$\quad almost everywhere
in $\Q$,
\end{enumerate}
where the operator ${\cal L}_0$ is defined by
\[ {\cal L}_0u:=u_t-\surm \div \Big[m(x) a_p (x,\n u)\Big]-\edeux f(u)=u_t-\Delta _\phi u -\edeux f(u).\]
In fact, if $\upl \in C^{2,1}(\Q)$ then, by Lemma
\ref{aniso-lapla-aniso}, the function ${\cal L}_0 \upl$ is
well-defined in $Q_T$. Also, using Lemma \ref{aniso-lapla-aniso},
we deduce that $\Delta _\phi  \upl \in L^\infty(Q_T)$. The
statement is then obtained by integrating
\eqref{deqexi-super-aniso} by parts. An analogous remark stands
for a sub-solution $\um\in C^{2,1}(\Q)$.\qed
\end{rem}

\vskip 8 pt We prove below an inequality which expresses the
strong monotonicity of the function $T^0(x,p)=a_p(x,p)$.
\begin{lem}
There exists a constant $\beta >0$ such that, for all
$x\in\ombar$, for all $p_1, p_2 \in \mathbb{R}^N$,
\begin{equation}\label{stromgmono-aniso}
(a_p(x,p_2)-a_p(x,p_1))\cdot(p_2-p_1)  \geq \beta  |p_2-p_1|^2.
\end{equation}
\end{lem}

{\noindent \bf Proof.} First we consider the case that $s
p_1+(1-s)p_2\neq 0$ for all $s\in [0,1]$. Then, the function
$s\mapsto a(x,s p_1+(1-s)p_2)$ is of class $C^2$ on $[0,1]$ and
there exist $p_3,p_4$ on the line segment $[p_1,p_2]$ such that
$$
a(x,p_2)-a(x,p_1)=a_p(x,p_1)\cdot(p_2-p_1)+\frac 1 2(p_2-p_1)\cdot
a_{pp}(x,p_3)(p_2-p_1),
$$
$$
a(x,p_1)-a(x,p_2)=a_p(x,p_2)\cdot(p_1-p_2)+\frac 1 2(p_1-p_2)\cdot
a_{pp}(x,p_4)(p_1-p_2).
$$
The strict convexity of $a(x,\cdot)$ implies that $a_{pp}(x,p)$ is
a positively definite symmetric matrix, so that the function
$(x,p,\bar p)\mapsto a_{pp}(x,p)\bar p \cdot \bar p$ is strictly
positive and continuous on the compact set $\ombar\times
S^{N-1}\times S^{N-1}$. Hence there exist constants $0<\lambda
_2\leq \Lambda _2$ such that, for all $x\in\ombar$, all $p\in
\mathbb{R}^N\setminus\{0\}$, all $\bar p\in \mathbb{R}^N$,
\begin{equation}\label{lambda2-aniso}
\lambda _2 |\bar p|^2\leq  a_{pp}(x,p)\bar p\cdot \bar p\leq
\Lambda _2 |\bar p|^2.
\end{equation}
It then follows that
\begin{align}
a(x,p_2)-a(x,p_1)\geq a_p(x,p_1)\cdot(p_2-p_1)+\frac {\lambda_2}
2|p_2-p_1|^2,\label{app1-aniso}\\
a(x,p_1)-a(x,p_2)\geq a_p(x,p_2)\cdot(p_1-p_2)+\frac {\lambda_2}
2|p_2-p_1|^2.\label{app2-aniso}
\end{align}
Adding up inequalities \eqref{app1-aniso} and \eqref{app2-aniso}
yields the desired inequality, with the constant $\beta=\lambda
_2$.

In the case that $sp_1+(1-s)p_2= 0$ for some $s\in [0,1]$, $p_1$
and $p_2$ are colinear and we may suppose that there exists $l\in
\R$ such that $p_2=lp_1$. We can assume $l\neq 0$, $l\neq 1$ and
$p_1 \neq 0$. By using the properties stated in Remark
\ref{der-a-aniso}, we obtain that
\begin{align*}
(a_p(x,p_2)-a_p(x,p_1))\cdot(p_2-p_1) &= (l-1)^2 a_p(x,p_1)\cdot
p_1\\
&= 2(l-1)^2 a(x,p_1)\\
&= 2 a(x,(l-1)p_1)\\
&\geq  {\lambda _0}^2 |(l-1)p_1|^2={\lambda _0}^2 |p_2-p_1|^2,
\end{align*}
where $\lambda _0$ has been defined in \eqref{lambda0-aniso}. The
proof is now completed. \qed

\vskip 8pt We are now ready to prove the following comparison
principle.
\begin{prop}[Comparison principle] Let $\upl$, respectively $\um$, be a super-solution,
 respectively a sub-solution, of Problem $\Pe$. Assume that
 $$
 \um(\cdot,0) \leq \upl(\cdot,0)\quad\text{almost everywhere in $\om$}.
 $$
Then we have that
\[ \um\leq \ue \leq \upl\quad\text{almost everywhere in $Q_T$.}\]
\end{prop}

{\noindent \bf Proof.} By subtracting inequality
\eqref{deqexi-super-aniso} for the super-solution $\upl$ from
inequality for the sub-solution $\um$, we obtain that, for all
$\varphi\in L^2(0,T;H^1(\om)) \cap L^\infty(Q_T)$ such that
$\varphi\geq 0$, and for all $t\in (0,T)$,
\begin{align}\label{dapp3-aniso}
\int_0^t\int_{\om}\Big[ (\um-\upl)_t\varphi+(a_p(x,\n
\um)&-a_p(x,\n \upl) )\cdot \n \varphi\Big]m(x)
\nonumber\\
&\leq  C\int_0^t\int_{\om}|\um- \upl|\varphi\,,
\end{align}
where $C$ is a constant depending on $\ep$ and the $L^\infty$
norms of $f'$ and $m$. Next we set $\varphi =(\um-\upl)^+$, which
belongs to $L^2(0,T;H^1(\om))\cap L^\infty(Q_T)$; it follows from
\eqref{stromgmono-aniso} that
\begin{align}
\int_0^t\int_{\om}&(a_p(x,\n \um)-a_p(x,\n \upl))\cdot \n
\varphi \,m(x)\nonumber \\
&=\int_0^t\int_{\{\um-\upl\geq 0\}}(a_p(x,\n \um)-a_p(x,\n
\upl))\cdot (\n\um-\n \upl)m(x)\nonumber \\
&\geq m_1\beta \int_0^t\int_{\{ \um-\upl\geq 0\}} |\n\um-\n
\upl|^2\geq 0\nonumber.
\end{align}
In view of (\ref{dapp3-aniso}), we now have that
$$
\di{\frac{m_1}{2}} \int_0^t \displaystyle{\frac{d}{dt}} \int_{\om}
\Big( (\um- \upl)^+\Big)^2 \leq C \int_0^t\int_{\{ \um-\upl\geq
0\}} (\um-\upl)^2\,,
$$
and therefore
\begin{equation*}
\int_{\om} \Big( (\um-\upl)^+\Big)^2(t)\leq \frac{2 C}{m_1}
\int_0^t\int_{\om} \Big( (\um-\upl)^+\Big)^2+\int_\om\Big(
(\um-\upl)^+\Big)^2(0).
\end{equation*}
Gronwall's lemma yields
\begin{equation*}
\int_{\om}\Big( (\um-\upl)^+\Big)^2(t)\leq
e^{2Ct/m_1}\int_\om\Big( (\um-\upl)^+\Big)^2(0).
\end{equation*}
Since $\um(x,0)\leq \upl(x,0)$ for almost all $x\in\om$, it
follows that
\[ \um \leq \upl \quad\text{a.e. in $Q_T$.}\] \qed

\begin{lem}
Let $\ue$ be the solution of Problem $\Pe$ (with initial data
$u_0$). Then
$$-\Vert u_0 \Vert_{L^\infty(\om)} \leq \ue \leq \max(1,\Vert u_0
\Vert_{L^\infty(\om)})\quad\  \text{a.e. in}\; \ Q_T.$$
\end{lem}

{\noindent \bf Proof.} By the bistable profile of $f$, we remark
that $-\Vert u_0 \Vert_{L^\infty(\om)}$, respectively
$\max(1,\Vert u_0 \Vert_{L^\infty(\om)})$, is a sub-solution,
respectively a super-solution, of Problem $\Pe$. \qed

\section{Generation of the interface}\label{generation-aniso}

This section deals with the generation of the interface, namely
the rapid formation of internal layers that takes place in a
neighborhood of $\Gamma_0=\{x\in \om,\,u_0(x)=a\}$ within the time
span of order $\ep^2 |\ln\ep|$. In the sequel, $\eta _0$ will
stand for the quantity
$$
\eta _0:= \min (a,1 -a).
$$
Our main result in this section is the following.

\begin{thm}\label{g-th-gen-aniso}
Let $\eta \in (0,\eta _0)$ be arbitrary and define $\mu$ as the
derivative of $f(u)$ at the unstable equilibrium $u=a$, that is
\begin{equation}\label{g-def-mu-aniso}
\mu=f'(a).
\end{equation}
Then there exist positive constants $\ep_0$ and $M_0$ such that,
for all $\,\ep \in (0,\ep _0)$,
\begin{enumerate}
\item for almost all $x\in\om$,
\begin{equation}\label{g-part1-aniso}
-\eta \leq u^\ep(x,\mu ^{-1}  \ep ^2 | \ln \ep |) \leq 1+\eta\,,
\end{equation}
\item for almost all $x\in\om$ such that $|u_0(x)-a|\geq M_0 \ep$,
we have that
\begin{align}
&\text{if}\;~~u_0(x)\geq a+M_0\ep\;~~\text{then}\;~~u^\ep(x,\mu
^{-1}  \ep ^2 | \ln \ep |)
\geq 1-\eta,\label{g-part2-aniso}\\
&\text{if}\;~~u_0(x)\leq a-M_0\ep\;~~\text{then}\;~~u^\ep(x,\mu
^{-1}  \ep ^2 | \ln \ep |)\leq \eta \label{g-part3-aniso}.
\end{align}
\end{enumerate}
\end{thm}

\vskip 8pt We will prove this result by constructing a suitable
pair of sub and super-solutions.

\subsection{The bistable ordinary differential equation}
The sub- and super-solutions mentioned above will be constructed
by modifying the solution of the problem without diffusion:
\begin{equation*}\label{no-diffusion-aniso}
\bar{u}_t=\frac{1}{\ep^2}\,f(\bar{u}), \qquad \bar{u}(x,0)=u_0(x).
\end{equation*}
This solution is written in the form
\[
\bar{u}(x,t)=Y\Big(\frac{t}{\ep^2},\,u_0(x)\Big),
\]
where $Y(\tau,\xi)$ denotes the solution of the ordinary
differential equation
\begin{equation}\label{g-ode-aniso}
\left\{\begin{array}{ll} Y_\tau (\tau,\xi)&=f(Y(\tau,\xi)) \quad\
\text{for} \; \ \tau >0 \vspace{3pt}\\
Y(0,\xi)&=\xi.
\end{array}\right.
\end{equation}
Here $\xi$ ranges over the interval $(-2C_0,2C_0)$, with $C_0$
being the constant defined in \eqref{int1-aniso}. We first collect
basic properties of $Y$.

\begin{lem}\label{g-Y1-aniso}
We have $Y_\xi >0$, for all $\xi \in (-2C_0,2C_0)$ and all $\tau >
0$. Furthermore,
$$
Y_\xi(\tau,\xi)=\frac{f (Y(\tau,\xi))}{f (\xi)}.
$$
\end{lem}

{\noindent \bf Proof.} First, differentiating equation
\eqref{g-ode-aniso} with respect to $\xi$, we obtain
\begin{equation}\label{nouveau1-aniso}
\left\{\begin{array}{ll} Y_{\xi\tau}=Y_\xi f'(Y),\vspace{3pt}\\
Y_\xi(0,\xi)=1,
\end{array}\right.
\end{equation}
which can be integrated as follows:
\begin{equation}\label{g-Y2-aniso}
Y_\xi(\tau,\xi)=\exp \Big[\int_0^\tau f'(Y(s,\xi))ds\Big]>0.
\end{equation}
We then differentiate equation \eqref{g-ode-aniso} with respect to
$\tau$ and obtain
\begin{equation}\label{nouveau2-aniso}
\left\{\begin{array}{ll} Y_{\tau\tau}=Y_\tau f'(Y),\vspace{3pt}\\
Y_\tau(0,\xi)=f (\xi),
\end{array}\right.
\end{equation}
which in turn implies
\begin{equation}\label{nouveau3-aniso}
Y_\tau(\tau,\xi)=f (\xi) \exp \Big[\di{\int_0^\tau}
f'(Y(s,\xi))ds\Big]=f (\xi)Y_\xi(\tau,\xi).
\end{equation}
This last equality, in view of \eqref{g-ode-aniso}, completes the
proof of Lemma \ref{g-Y1-aniso}. \qed

\vskip 8pt We define a function $A(\tau,\xi)$ by
\begin{equation}\label{g-A-aniso}
A(\tau,\xi)=\frac{f'(Y(\tau,\xi))-f'(\xi)}{f (\xi)}.
\end{equation}

\begin{lem}\label{g-Y5-aniso}
We have, for all $\xi \in (-2C_0,2C_0)$ and all $\tau > 0$,
$$
A(\tau,\xi)=\int_0^\tau f''(Y(s,\xi))Y_\xi(s,\xi)ds.
$$
\end{lem}

{\noindent \bf Proof.} Differentiating the equality of Lemma
\ref{g-Y1-aniso} with respect to $\xi$  leads to
\begin{equation}\label{g-Y4-aniso}
Y_{\xi\xi}=A(\tau,\xi)Y_\xi,
\end{equation}
whereas differentiating \eqref{g-Y2-aniso} with respect to $\xi$
yields
$$
Y_{\xi\xi}=Y_\xi \int_0^\tau f''(Y(s,\xi))Y_\xi(s,\xi)ds.
$$
These two last results complete the proof of Lemma
\ref{g-Y5-aniso}. \qed

\vskip 8 pt Next we need some estimates on $Y$ and its
derivatives. First, we perform some estimates when the initial
value $\xi$ lies between $\eta$ and $1-\eta$.

\vskip 8pt
\begin{lem}\label{g-est-derY-A-milieu-aniso}
Let $\eta \in (0,\eta _0)$ be arbitrary. Then there exist positive
constants $\tilde C _1=\tilde C _1(\eta)$, $\tilde C _2=\tilde C
_2(\eta)$ and $C_3=C_3(\eta)$ such that, for all $\tau>0$,
\begin{enumerate}
\item if $\xi \in (a,1-\eta)$ then, for every $\tau >0$ such that
$Y(\tau,\xi)$ remains in the interval $(a,1-\eta)$, we have
\begin{equation}\label{g-est-Y3-aniso}
\tilde C _1e^{\mu\tau}\leq Y_\xi(\tau,\xi) \leq \tilde C _2
e^{\mu\tau}\,,
\end{equation}
and
\begin{equation}\label{g-est-A-milieu-aniso}
|A(\tau,\xi)|\leq C_3(e^{\mu \tau}-1);
\end{equation}
\item if $\xi\in (\eta,a)$ then, for every $\tau >0$ such that
$Y(\tau,\xi)$ remains in the interval $(\eta,a)$,
\eqref{g-est-Y3-aniso} and \eqref{g-est-A-milieu-aniso} hold as
well,
\end{enumerate}
where $\mu$ is the constant defined in \eqref{g-def-mu-aniso}.
\end{lem}

{\noindent \bf Proof.} We take $\xi \in (a,1-\eta)$ and suppose
that for $s \in (0,\tau)$, $Y(s,\xi)$ remains in the interval
$(a,1-\eta)$. Integrating the equality $Y_\tau / f(Y)=1$ from $0$
to $\tau$ yields
\begin{equation}\label{g-tau-aniso}
\tau=\int_0^{\tau} \frac{Y_\tau (s,\xi)}{f (Y(s,\xi))}ds=\int _\xi
^{Y(\tau,\xi)} \frac{dq}{f (q)}.
\end{equation}
Moreover, the equality of Lemma \ref{g-Y1-aniso} leads to
\begin{equation}\label{naka-aniso}
\begin{array}{lll}
\ln Y_\xi (\tau,\xi)&=
\di{\int _ \xi ^{Y(\tau,\xi)}} \frac{f'(q)}{f(q)}dq \vsp \\
&=\di{\int _ \xi ^{Y(\tau,\xi)}}\big
[\frac{f'(a)}{f(q)}+\frac{f'(q)-f'(a)}{f(q)}\big ]dq \vsp \\
&=\mu \tau+\di{\int _ \xi ^{Y(\tau,\xi)}}h(q)dq,
\end{array}
\end{equation}
where
$$
h(q)=(f'(q)-\mu)/f(q).
$$
Since
$$
h(q)\to\frac{f''(a)}{f'(a)} \quad\ \text{as}\; \ q\to a,
$$
the function $h$ is continuous on $[a,1-\eta]$. Hence we can
define
$$
H=H(\eta):=\Vert h \Vert _{L^\infty (a,1-\eta)}.
$$
Since $|Y(\tau,\xi)-\xi|$ takes its values in the interval $[0,1
-a-\eta]\subset[0,1-a]$, it follows from \eqref{naka-aniso} that
$$
\mu \tau -H(1-a) \leq \ln Y_\xi(\tau,\xi) \leq \mu \tau+H(1-a),
$$
which, in turn, proves \eqref{g-est-Y3-aniso}. Lemma
\ref{g-Y5-aniso} and \eqref{g-est-Y3-aniso} yield
$$
|A(\tau,\xi)| \leq \sup_{z \in [0,1]} |f'' (z)| \di{\int_0^\tau}
\tilde C _2 e^{\mu s}ds\leq C_3(e^{\mu\tau}-1),
$$
which completes the proof of \eqref{g-est-A-milieu-aniso}. The
case where $\xi$ and $Y(\tau,\xi)$ are in $(\eta,a)$ is similar
and omitted. \qed

\begin{cor}\label{g-est-Y-milieu-aniso} Let $\eta \in (0,\eta
_0)$ be arbitrary. Then there exist positive constants
$C_1=C_1(\eta)$ and $C_2=C_2(\eta)$ such that, for all $\tau>0$,
\begin{enumerate}
\item if $\xi\in (a,1-\eta)$ then, for every $\tau >0$ such that
$Y(\tau,\xi)$ remains in the interval $(a,1-\eta)$, we have
\begin{equation}\label{g-est-Y-1-aniso}
C_1e^{\mu \tau}(\xi-a)\leq Y(\tau,\xi)-a \leq C_2e^{\mu
\tau}(\xi-a);
\end{equation}
\item if $\xi \in (\eta,a)$ then, for every $\tau >0$ such that
$Y(\tau,\xi)$ remains in the interval $(\eta,a)$, we have
\begin{equation}\label{g-est-Y-2-aniso}
C_2e^{\mu \tau}(\xi-a)\leq Y(\tau,\xi)-a \leq C_1e^{\mu
\tau}(\xi-a).
\end{equation}
\end{enumerate}
\end{cor}

{\noindent \bf Proof.} Since
$$
f(q)/(q-a)\to f'(a) \quad\ \text{as}\; \ q \to a,
$$
it is possible to find $B_1=B_1(\eta)>0$ and $B_2=B_2(\eta)>0$
such that, for all $q\in (a,1-\eta)$,
\begin{equation}\label{g-ineg-fq-aniso}
B_1(q-a)\leq f (q) \leq B_2(q-a).
\end{equation}
We write this inequality for $a<Y(\tau,\xi) < 1-\eta$ to obtain
$$
B_1(Y(\tau,\xi)-a)\leq f (Y(\tau,\xi))\leq B_2(Y(\tau,\xi)-a).
$$
We also write this inequality for $a<\xi < 1-\eta$ to obtain
$$
B_1(\xi-a) \leq f(\xi) \leq B_2(\xi-a).
$$
Next we use the equality $Y_\xi=f (Y)/f (\xi)$ of Lemma
\ref{g-Y1-aniso} to deduce that
$$
\frac{B_1}{B_2}(Y(\tau,\xi)-a)\leq (\xi-a)Y_\xi(\tau,\xi) \leq
\frac{B_2}{B_1}(Y(\tau,\xi)-a),
$$
which, in view of \eqref{g-est-Y3-aniso}, implies that
$$
\frac{B_1}{B_2}\tilde C_1 e^{\mu \tau}(\xi-a) \leq Y(\tau,\xi)-a
\leq \frac{B_2}{B_1}\tilde C_2 e^{\mu \tau}(\xi-a).
$$
This proves \eqref{g-est-Y-1-aniso}. The proof of
\eqref{g-est-Y-2-aniso} is similar and omitted. \qed

\vskip 8 pt Next we present estimates in the case where the
initial value $\xi$ is smaller than $\eta$ or larger than
$1-\eta$.

\begin{lem}\label{g-est-bords-aniso}
Let $\eta \in (0,\eta _0)$ and $M>0$ be arbitrary.  Then there
exists a positive constant $C_4=C_4(\eta,M)$ such that
\begin{enumerate}
\item if $\xi \in [1-\eta,1+M]$, then, for all $\tau> 0$,
$Y(\tau,\xi)$ remains in the interval $[1-\eta,1+M]$ and
\begin{equation}\label{g-est-A-bords-aniso}
|A(\tau,\xi)|\leq C_4\tau \quad\hbox{for}\ \ \tau> 0 \,;
\end{equation}
\item if $\xi \in [-M,\eta]$, then, for all $\tau> 0$,
$Y(\tau,\xi)$ remains in the interval $[-M,\eta]$ and
\eqref{g-est-A-bords-aniso} holds as well.
\end{enumerate}
\end{lem}

{\noindent \bf Proof.} Since the two statements can be treated in
the same way, we will only prove the former. The fact that
$Y(\tau,\xi)$, the solution of the ordinary differential equation
\eqref{g-ode-aniso}, remains in the interval $[1-\eta,1+M]$
directly follows from the bistable properties of $f$, or, more
precisely, from the sign conditions $f(1-\eta)>0$, $f(1+M)<0$.

To prove \eqref{g-est-A-bords-aniso}, suppose first that $\xi\in
[1,1+M]$. In view of \eqref{der-f-aniso}, $f'$ is strictly
negative in an interval of the form $[1,1 +c]$ and $f$ is negative
in $[1,\infty)$.  We denote by $-m<0$ the maximum of $f$ on
$[1+c,1+M]$.  Then, as long as $Y(\tau,\xi)$ remains in the
interval $[1+c,1+M]$, the ordinary differential equation
\eqref{g-ode-aniso} implies
$$
Y_\tau \leq -  m .
$$
By integration, this means that, for any $\xi\in [1,1+M]$, we have
\[
Y(\tau,\xi) \in [1,1 +c]\qquad\hbox{for}\ \ \tau
\geq\overline{\tau}:=\frac{M-c}{m}.
\]
In view of this, and considering that $f'(Y)<0$ for $Y\in
[1,1+c]$, we see from the expression \eqref{g-Y2-aniso} that
\[
\begin{array}{ll}
Y_\xi(\tau,\xi) &\di\vspace{6pt}
=\exp\Big[\int_0^{\overline{\tau}}f'(Y(s,\xi))ds\Big]\,
\exp\Big[\int_{\overline{\tau}}^\tau f'(Y(s,\xi))ds\Big]\\
&\di\vspace{6pt}
\leq \exp \Big[\int _0^{\overline{\tau}}f'(Y(s,\xi))ds\Big]\\
&\di \leq \exp\Big[\int _0 ^{\overline{\tau}} \sup_{z\in[-M,1+M]}
|f'(z)|ds\Big] =:\tilde{C}_4=\tilde{C}_4(M),
\end{array}
\]
for all $\tau\geq\overline{\tau}$. It is clear from the same
expression \eqref{g-Y2-aniso} that $Y_\xi\leq\tilde{C}_4$ holds
also for $0\leq\tau\leq\overline{\tau}$.  We can then use Lemma
\ref{g-Y5-aniso} to deduce that
$$
\begin{array}{ll}
|A(\tau,\xi)|&\leq \tilde{C}_4 \di{\int_0^\tau} |f''(Y(s,\xi))|ds\vsp\\
&\leq \tilde{C}_4 \Big(\sup _{z\in[-M,1+M]} |f''(z)|\Big)
\tau=:C_4\tau.
\end{array}
$$
The case $\xi\in [1-\eta,1]$ can be treated in the same way. This
completes the proof of the lemma. \qed

\vskip 8 pt Now we choose the constant $M$ in the above lemma
sufficiently large so that $[-2C_0,2C_0]\subset [-M,1+M]$, and fix
$M$ hereafter. Then $C_4$ only depends on $\eta$. Using the fact
that $\tau=O(e^{\mu\tau}-1)$ for $\tau>0$, one can easily deduce
from \eqref{g-est-A-milieu-aniso} and \eqref{g-est-A-bords-aniso}
the following general estimate.

\begin{lem}\label{g-EST-A-aniso}
Let $\eta \in (0,\eta _0)$ be arbitrary and let $C_0$ be the
constant defined in \eqref{int1-aniso}. Then there exists a
positive constant $C_5=C_5(\eta)$ such that, for all $\xi
\in(-2C_0,2C_0)$ and all $\tau>0$,
$$
|A(\tau,\xi)|\leq C_5(e^{\mu \tau}-1).
$$
\end{lem}

\subsection{Construction of sub- and super-solutions}

We are now ready to construct sub- and super-solutions in order to
study the generation of the interface. By using some cut-off
initial data, see subsection 3.2 in \cite{AHM}, we can modify
slightly $u_0$ near the boundary $\partial\Omega$ and make,
without loss of generality, the additional assumption
\begin{equation}\label{int2-aniso}
a_p(x,\n u_0(x))\cdot\nu =0 \quad\ \text{on}\; \ \partial\Omega.
\end{equation}
Our sub- and super-solutions are defined by
\begin{equation}\label{w+--aniso}
w_\ep^\pm(x,t)=
Y\Big(\frac{t}{\ep^2},\,u_0(x)\pm\ep^2C_6(\emuth-1)\Big).
\end{equation}

\begin{lem}\label{g-w-aniso}
There exist positive constants $\ep_0$ and $C_6$ such that, for
all $\, \ep \in (0,\ep _0)$, $(w_\ep^-,w_\ep^+)$ is a pair of sub-
and super-solutions of Problem $\Pe$, in the domain $\om \times
(0,\mu ^{-1} \ep^2|\ln \ep|)$.
\end{lem}

{\noindent \bf Proof.} Following Remark \ref{rem-sup-aniso} we
define the operator ${\cal L}_0$ by
\begin{equation}\label{operator-aniso}
{\cal L}_0 u :=u_t -\surm \div\Big [ m(x) a_p(x,\n u)\Big] -\edeux
f(u)\,,
\end{equation}
and prove that ${\cal L}_0 w_\ep^+\geq 0$. We compute
\begin{align*}
({w_\ep}^+)_t&=\di{\frac{1}{\ep^2}}Y_\tau+\mu C_6 \emuth Y_\xi,\\
\n w_ \ep ^+&=\n u_0(x) Y_\xi.
\end{align*}
Using  \eqref{int2-aniso} and the fact that $a_p(x,\cdot)$ is
homogeneous of degree one, we see that $w_\ep ^\pm$ satisfy the
anisotropic Neumann boundary condition $ a_p(x,\n w_\ep^\pm)\cdot
\nu=0$ on $\partial \Omega \times (0,+\infty)$. In view of the
ordinary differential equation \eqref{g-ode-aniso}, we obtain
$$
{\cal L}_0 w_\ep^+ =\mu C_6 \emuth Y_\xi -\surm \div \Big [ m(x)
a_p(x,\n w_\ep ^+)\Big].
$$
By the estimate of the anisotropic Laplacian
\eqref{aniso-lapla-eq-aniso}, it follows that
\begin{equation}\label{etape-aniso}
{\cal L}_0w_\ep^+\geq \mu C_6 \emuth Y_\xi-C_L(|\n w_\ep
^+(x,t)|+|D^2 w_\ep ^+(x,t)|)\,,
\end{equation}
where we recall that $|D^2 w_\ep ^+(x,t)|=\max_{i,j}|\partial _i
\partial _j w _\ep ^+ (x,t)|$. A straightforward calculation yields
$$
\partial _i\partial _j w _\ep ^+ (x,t)=(\partial _i\partial _j u_0)Y_\xi+(\partial _i u_0\partial _j
u_0)Y_{\xi\xi}.
$$
Recalling that $Y_\xi>0$, we now combine the expression of $\n
w_\ep^+$, the above expression and inequality \eqref{etape-aniso}
to obtain
\begin{align}\label{intermerdiaire-aniso}
{\cal L}_0w_\ep^+ /Y_\xi \geq \mu
C_6\emuth-C_LC_0-C_0-{C_0}^2\di{\frac{|Y_{\xi\xi}|}{Y_\xi}}\,,
\end{align}
where $C_0$ is the constant defined in \eqref{int1-aniso}. We note
that, in the range $(0,\mu ^{-1} \ep ^2|\ln \ep|)$, we have
$$
0\leq \ep^2 C_6(\emuth-1) \leq \ep ^2C_6(\ep^{-1}-1) \leq C_0\,,
$$
if $\ep _0$ is small enough. Hence
$$
\xi:=u_0(x)\pm \ep ^2 \ep C_6(\emuth-1)\in (-2C_0,2C_0),
$$
so that, by the results of the previous subsection, $Y$ remains in
$(-2C_0,2C_0)$. In view of \eqref{g-Y4-aniso}, $Y_{\xi\xi}/Y_\xi$
is equal to $A$ so that, combining the estimate of $A$ in Lemma
\ref{g-EST-A-aniso} and \eqref{intermerdiaire-aniso}, we obtain
$$
{\cal L}_0 w_\ep^+/Y_\xi \geq ( \mu C_6
-{C_0}^2C_5)\emuth-C_LC_0-C_0.
$$
Now, choosing
$$
C_6\geq\frac 2 \mu \max\big({C_0}^2C_5,C_0(C_L+1)\big)
$$
proves ${\cal L}_0 w_\ep^+/Y_\xi\geq 0$. Since $Y_\xi>0$, it
follows that ${\cal L}_0 w_\ep^+\geq 0$. Hence, by Remark
\ref{rem-sup-aniso}, $w_\ep^+$ is a super-solution of Problem
$\Pe$. Similarly $w_\ep^-$ is a sub-solution. Lemma
\ref{g-w-aniso} is proved. \qed

\vskip 8 pt To conclude this subsection, we remark that
$w^\pm(x,0)=Y\Big(\di{\frac{t}{\ep^2}},\,u_0(x)\Big)=u_0(x)$.
Consequently, by the comparison principle,
\begin{equation}\label{g-coincee1-aniso}
w_\ep^-(x,t) \leq u^\ep(x,t) \leq w_\ep^+(x,t)\,,
\end{equation}
for almost all $(x,t)\in \om\times(0,\mu ^{-1} \ep ^2 |\ln \ep|)$.

\subsection{Proof of Theorem \ref{g-th-gen-aniso}}

In order to prove Theorem \ref{g-th-gen-aniso} we first present a
key estimate on the function $Y$ after a time interval of order
$\tau\sim |\ln \ep|.$

\begin{lem}\label{after-time-aniso}
Let $\eta \in (0,\eta _0)$ be arbitrary; there exist positive
constants $\ep_0$ and $C_7$ such that, for all $\,\ep\in(0,\ep
_0)$,
\begin{enumerate}
\item for all $\xi\in (-2C_0,2C_0)$, for all $\tau \geq \mu
^{-1}|\ln \ep|$,
\begin{equation}\label{g-part11-aniso}
-\eta \leq Y(\tau,\xi) \leq 1+\eta\,,
\end{equation}
\item for all $\xi\in (-2C_0,2C_0)$ such that $|\xi-a|\geq C_7
\ep$, for all $\tau \geq \mu ^{-1}|\ln \ep|$,
\begin{align}
&\text{if}\;~~\xi\geq a+C_7 \ep\;~~\text{then}\;~~Y(\tau,\xi)
\geq 1-\eta,\label{g-part22-aniso}\\
&\text{if}\;~~\xi\leq a-C_7 \ep\;~~\text{then}\;~~Y(\tau,\xi)\leq
\eta \label{g-part33-aniso}.
\end{align}
\end{enumerate}
\end{lem}

{\noindent \bf Proof.} We first prove \eqref{g-part22-aniso}. For
$\xi \geq a+C_7\ep$, as long as $Y(\tau,\xi)$ has not reached
$1-\eta$, we can use \eqref{g-est-Y-1-aniso} to deduce that
$$
Y(\tau,\xi) \geq a +C_1 e^{\mu \tau} (\xi-a) \geq a + C_1 C_7
e^{\mu \tau} \ep \geq 1-\eta\,,
$$
provided that $\tau$ satisfies $\tau \geq \mu ^{-1}\ln \frac {1
-a-\eta}{C_1C_7 \ep}$. Choosing
$$
C_7=\frac{\max (a,1-a)-\eta}{C_1}
$$
completes the proof of \eqref{g-part22-aniso}. Using
\eqref{g-est-Y-2-aniso}, one easily proves \eqref{g-part33-aniso}.

Next we prove \eqref{g-part11-aniso}. First, by the bistable
assumptions on $f$, if we leave from a $\xi \in [-\eta,1+\eta]$
then $Y(\tau,\xi)$ will remain in $[-\eta,1+\eta]$. Now suppose
that $1+\eta \leq \xi \leq 2C_0$. We check below that $Y(\mu
^{-1}|\ln \ep|,\xi)\leq 1+\eta$. First, in view of
\eqref{der-f-aniso}, we can find $p>0$ such that
\begin{equation}\label{g-pente-aniso}
\begin{array}{ll}\text { if } \quad 1 \leq
u \leq 2C_0 & \text { then } \quad  f(u)
\leq p(1-u),\vspace{3pt}\\
\text { if } \quad -2C_0 \leq u \leq 0 & \text { then } \quad f(u)
\geq -pu.
\end{array}
\end{equation}
We then use the ordinary differential equation \eqref{g-ode-aniso}
to obtain, as long as $1+\eta \leq Y \leq 2C_0$, the inequality
$Y_\tau \leq p(1-Y)$. It follows that
$$
\frac{Y_\tau}{Y-1}\leq -p.
$$
Integrating this inequality from $0$ to $\tau$ leads to
$$
Y(\tau,\xi)\leq 1+(\xi-1)e^{-p\tau}\leq 1+(2C_0-1)e^{-p\tau}.
$$
One easily checks that, for $\ep \in (0,\ep _0)$, with $\ep_0=\ep
_0(\eta)$ small enough, we have $Y(\tau,\xi)\leq 1+\eta$, for all
$\tau \geq \mu ^{-1}| \ln \ep|$, which completes the proof of
\eqref{g-part11-aniso}. \qed

\vskip 8 pt We are now ready to prove Theorem
\ref{g-th-gen-aniso}. By setting $t=\mu ^{-1} \ep ^2|\ln \ep|$ in
\eqref{g-coincee1-aniso}, we obtain, for almost all $x\in\om$,
\begin{align}
Y(\mu ^{-1}&|\ln \ep|, u_0(x)-(C_6 \ep -C_6 \ep ^2
))\nonumber\vsp\\
&\leq u^\ep(x,\mu ^{-1} \ep^2|\ln \ep|) \leq Y(\mu ^{-1}|\ln \ep|,
u_0(x)+C_6 \ep -C_6 \ep ^2).\label{g-gr-aniso}
\end{align}
Furthermore, by the definition of $C_0$ in \eqref{int1-aniso}, we
have, for $\ep_0$ small enough,
$$
- 2 C_0 \leq u_0(x)\pm (C_6 \ep -C_6 \ep ^2) \leq 2C_0
\qquad\hbox{for}\ x\in\om.
$$
Thus the assertion \eqref{g-part1-aniso} of Theorem
\ref{g-th-gen-aniso} is a direct consequence of
\eqref{g-part11-aniso} and \eqref{g-gr-aniso}.

Next we prove \eqref{g-part2-aniso}. We choose $M_0$ large enough
so that $M_0\ep- C_6 \ep+C_6 \ep ^2 \geq C_7\ep$.  Then, for any
$x\in \om$ such that $u_0(x)\geq a+M_0 \ep$, we have
$$
u_0(x)-(C_6 \ep -C_6 \ep ^2)\geq a+M_0\ep- C_6 \ep+C_6 \ep ^2 \geq
a+C_7\ep.
$$
Combining this, \eqref{g-gr-aniso} and \eqref{g-part22-aniso}, we
see that
\[
u^\ep(x,\mu^{-1}\ep ^2 | \ln \ep |)\geq 1-\eta\,,
\]
for almost all $x\in\om$ that satisfies $u_0(x)\geq a+M_0 \ep$.
This proves \eqref{g-part2-aniso}. The inequality
\eqref{g-part3-aniso} can be shown the same way.  This completes
the proof of Theorem \ref{g-th-gen-aniso}.\qed

\section{Motion of the interface}\label{motion-aniso}

We have seen in Section \ref{generation-aniso} that, after a very
short time, the solution $\ue$ develops a clear transition layer.
In the present section, we show that it persists and that its law
of motion is well approximated by the interface equation $(P^0)$.

More precisely, take the first term of the formal asymptotic
expansion \eqref{inner-aniso} as a formal expansion of the
solution:
\begin{equation}\label{approxmotion-aniso}
\ue(x,t)\,\approx\,\tilde{u}^\ep(x,t):= \U\left(\frac{\dtip
(x,t)}{\ep}\right).
\end{equation}
The right-hand side of \eqref{approxmotion-aniso} is a function
having a well-developed transition layer, and its interface lies
exactly on $\Gamma _t$. We show that this function is a very good
approximation of the solution; therefore the following holds:
\begin{quote}
{\it If $u^\ep$ becomes rather close to $\tilde{u}^\ep$ at some
time moment, then it stays close to $\tilde{u}^\ep$ for the rest
of time.}
\end{quote}

To that purpose, we will construct a pair of sub- and
super-solutions $u_\ep^-$ and $u_\ep^+$ of Problem $(P^\ep)$ by
slightly modifying $\tilde{u}^\ep$. It then follows that, if the
solution $\ue$ satisfies
\[
u_\ep^-(x,t_0)\leq \ue(x,t_0)\leq  u_\ep^+(x,t_0)\,,
\]
for some $t_0\geq 0$ and for almost all $x\in\om$, then
\[
u_\ep^-(x,t)\leq \ue(x,t)\leq  u_\ep^+(x,t)\,,
\]
for almost $(x,t)\in Q_T$ that satisfies $t_0\leq t\leq T$. As a
result, since both $u_\ep^+, u_\ep^-$ stay close to
$\tilde{u}^\ep$, the solution $\ue$ also stays close to
$\tilde{u}^\ep$ for $t_0\leq t\leq T$.

\subsection{Construction of sub and super-solutions}
To begin with we present a mathematical tool which is essential
for the construction of sub and super-solutions.

\vskip 8 pt {\noindent \bf A modified anisotropic signed distance
function.} Rather than working with the anisotro\-pic signed
distance function $\widetilde d _\phi$, defined in
\eqref{eq:dist-aniso}, we define a \lq \lq cut-off anisotropic
signed distance function" $d_ \phi$ as follows. Choose $d_0>0$
small enough so that $\widetilde d _\phi (\cdot,\cdot)$ is smooth
in the tubular neighborhood of $\Gamma$
$$
\{(x,t) \in {\Q},\;|\widetilde d _\phi(x,t)|<3d_0\},
$$
and that
\begin{equation}\label{front-aniso}
\mbox {dist} _\phi(\Gamma_t,\partial \Omega)>3d_0 \quad \textrm{
for all } t\in (0,T).
\end{equation}
Next let $\zeta(s)$ be a smooth increasing function on $\R$ such
that
\[
\zeta(s)= \left\{\begin{array}{ll}
s &\textrm{ if }\ |s| \leq d_0\vspace{4pt}\\
-2d_0 &\textrm{ if } \ s \leq -2d_0\vspace{4pt}\\
2d_0 &\textrm{ if } \ s \geq 2d_0.
\end{array}
\right.
\]
We define the cut-off anisotropic signed distance function $d
_\phi$ by
\begin{equation}
d _\phi(x,t)=\zeta\big(\widetilde d_\phi(x,t)\big).
\end{equation}
Note that, in view of \eqref{bellettini-aniso},
\begin{equation}\label{bellettini-motion-aniso}
2a(x,\n d_\phi(x,t))=1 \quad \quad \text{ in a neighborhood of }
\Gamma _t,
\end{equation}
more precisely in the region $\{(x,t) \in {\Q},\,|d
_\phi(x,t)|<d_0\}$. Moreover, in view of \eqref{front-aniso}, we
have
\begin{equation}\label{pres-bord-aniso}
2a(x,\n d_\phi(x,t))=0 \quad\quad \text{ far away from } \Gamma
_t,
\end{equation}
more precisely in the region $\{(x,t) \in {\Q},\,|d
_\phi(x,t)|\geq 2d_0\}$. Furthermore, since the moving interface
$\Gamma$ satisfies Problem $(P^0)$, an alternative equation for
$\Gamma$ is given by
\begin{equation}\label{FBP-aniso}
(d _\phi) _t=\surm \div \big[m(x) a_p (x,\n d _\phi)\big] \quad
\text{ on }\Gamma _t.
\end{equation}

\vskip 8 pt {\noindent \bf Construction.} We look for a pair of
sub- and super-solutions $u_\ep^{\pm}$ for $\Pe$ of the form
\begin{equation}\label{sub-aniso}
u_\ep^{\pm}(x,t)=U_0\Big(\frac{d _\phi(x,t) \pm \ep
p(t)}{\ep}\Big) \pm q(t),
\end{equation}
where $U_0$ is the solution of \eqref{eq-phi-aniso}, and where
$$
\begin{array}{ll}
p(t)=-\EB+e^{Lt}+ K,\vsp \\
q(t)=\sigma(\beta \EB+\ep^2Le^{Lt}).
\end{array}
$$
Note that $q=\sigma\ep^2\,p_t$.  It is clear from the definition
of $u_\ep^\pm$ that
\begin{equation}\label{sub2-aniso}
\lim_{\ep\rightarrow 0} u_\ep^\pm(x,t)= \left\{
\begin{array}{ll}
1 &\textrm { for all } (x,t) \in Q_T^+ \vspace{4pt}\\
0 &\textrm { for all } (x,t) \in Q_T^-.\\
\end{array}\right.
\end{equation}

The main result of this section is the following.
\begin{lem}\label{fix-aniso}
There exist positive constants $\beta,\,\sigma$ with the following
properties.  For any $K>1$, we can find positive constants $\ep_0$
and $L$ such that, for any $\ep\in(0,\ep _0)$, the functions
$u_\ep^-$ and $u_\ep^+$ satisfy the anisotropic Neumann boundary
condition and
$$
{\cal L}_0u_\ep^-\leq 0 \leq {\cal L}_0u_\ep^+,
$$
in the range $\om\times(0,T)$, where the operator ${\cal L}_0$ has
been defined in \eqref{operator-aniso}.
\end{lem}

\subsection{Proof of Lemma \ref{fix-aniso}}

We show below that $${\cal L}_0 u_\ep^+ :=(u_\ep^+)_t -\surm \div
\big[m(x) a_p(x,\n u_\ep^+)\big] -\edeux f(u_\ep^+)\geq 0,$$ the
proof of inequality ${\cal L}_0 u_\ep^-\leq 0$ follows by similar
arguments.

\subsubsection {Computation of ${\cal L}_0 u_\ep^+$}
In the sequel, the function $U_0$ and its derivatives are taken at
the point $(d _\phi (x,t)+\ep p(t))/ \ep$. Straightforward
computations yield
\begin{align*}
(u_\ep^+)_t&= (\frac{1}{\ep}(d _\phi) _t+p_t){U_0}'+q_t,\\
\nabla u_\ep^+ &= \frac{1}{\ep}{U_0}'\nabla d _\phi,\\
\div a_p(x,\n u_\ep ^+)&=\edeux {U_0}'' \n \dphi \cdot a_p(x,\n
\dphi)+\frac 1 \ep {U_0}'\div a_p(x,\n \dphi)\\
&=\edeux {U_0}'' 2a(x,\n \dphi)+\frac 1 \ep {U_0}'\div a_p(x,\n
\dphi),
\end{align*}
where we have used properties stated in Remark \ref{der-a-aniso}.
Note that, $d _\phi$ being constant in a neighborhood of $\partial
\Omega $, we have that $\n u_\ep^+=0$ on $\partial\Omega
\times(0,T)$ and $u _\ep ^+$ satisfies the anisotropic Neumann
boundary condition $a_p (x,\n u_\ep ^+)\cdot\nu =0$ on $\partial
\om \times (0,T)$. At last, we use the expansion
$$
f(u_\ep^+)=f(U_0)+qf'(U_0)+\frac 12 q^2f''(\theta)\,,
$$
for some function $\theta(x,t)$ satisfying $U_0<\theta<u_\ep^+$.

Combining the above expressions with \eqref{mdiv} and
\eqref{eq-phi-aniso-pb}, we obtain ${\cal L}_0
u_\ep^+=E_1+E_2+E_3$, where
\begin{align*}
E_1&=- \edeux  q\big(f'(U_0)+\frac 12 q f''(\theta)\big)+{U_0}' p_t+q_t,\vsp\\
E_2&=\frac{{U_0}''}{\ep^2}\Big(1-2a(x,\n \dphi)\Big),\vsp\\
E_3&=\frac {{U_0}'}{\ep}\Big((d _\phi) _t-\surm \div \big[m(x)
a_p(x,\n \dphi)\big]\Big).
\end{align*}

\vskip 8 pt In order to estimate the above terms, we first present
some useful inequalities. As $f'(0)$ and $f'(1)$ are strictly
negative, we can find strictly positive constants $b$ and $m$ such
that
\begin{equation}\label{bords-aniso}
\textrm{ if }\quad U_0(z) \in [0,b]\cup [1-b,1] \quad \quad
\textrm{ then }\quad f'(U_0(z))\leq -m.
\end{equation}
On the other hand, since the region $\{(x,z)\in\ombar \times
\R\,|\,U_0(z)\in [b,1-b] \,\}$ is compact and since ${U_0}'>0$ on
$\R$, there exists a constant $a_1>0$ such that
\begin{equation}\label{milieu-aniso}
\textrm{ if }\quad U_0(z)\in [b,1-b] \quad \textrm{ then } \quad
{U_0}'(z)\geq a_1.
\end{equation}
We now choose $M>0$ such that $|U_0|\leq M -1$. We then define
\begin{equation}\label{F-aniso}
F=\sup_{|z|\leq M} |f(z)|+|f'(z)|+|f ''(z)|\, ,
\end{equation}
\begin{equation}\label{beta-aniso}
\beta = \frac{m }{4}\,,
\end{equation}
and choose $\sigma$ that satisfies
\begin{equation}\label{sigma-aniso}
0<\sigma\leq\min(\sigma _0, \sigma _1, \sigma _2)\,,
\end{equation}
where
$$
\sigma _0:=\frac{a_1}{4\beta+F}, \quad \sigma _1:=\frac{1}{\beta
+1}, \quad \sigma _2:=\frac{4\beta }{F(\beta+1)}.
$$
Hence, combining \eqref{bords-aniso} and \eqref{milieu-aniso}, we
obtain, using that $\sigma \leq \sigma _0$,
\begin{equation}\label{U0-f-aniso}
{U_0}'(z)-\sigma f'(U_0(z))\geq 4 \sigma \beta \qquad
\hbox{for}\,\,\; z\in\R.
\end{equation}

Now let $K>1$ be arbitrary. In what follows we will show that
${\cal L}_0 u_\ep^+\geq 0$ provided that the constants $\ep_0$ and
$L$ are appropriately chosen. From now on, we suppose that the
following inequality is satisfied:
\begin{equation}\label{ep0M-aniso}
\ep _0^2 Le^{LT} \leq 1\, .
\end{equation}
Then, given any $\ep\in(0,\ep_0)$, since $\sigma \leq \sigma _1$,
we have $0\leq q(t)\leq 1$, hence
\begin{equation}\label{uep-pm-aniso}
-M\leq u_\ep^\pm(x,t) \leq M\, .
\end{equation}

\subsubsection {An estimate for $E_1$}

A direct computation gives
$$
E_1=\frac{\beta}{\ep^2}\,\EB(I-\sigma\beta)+Le^{Lt}(I+\ep^2\sigma
L),
$$
where
$$
I={U_0}'-\sigma  f '(U_0)-\frac {\sigma^2}2  f
''(\theta)(\beta\EB+\ep^2 Le^{Lt}).
$$
In virtue of \eqref{U0-f-aniso} and \eqref{uep-pm-aniso}, we
obtain
\[
I\geq 4 \sigma \beta-\frac {\sigma^2}{2}F(\beta+\ep^2 Le^{LT}).
\]
Then, in view of \eqref{ep0M-aniso}, using that $\sigma \leq
\sigma _2$, we have $I\geq 2\sigma\beta$, which implies
$$
E_1\geq \frac{\sigma\beta^2}{\ep^2}\EB + 2\sigma\beta L
e^{Lt}=:\frac {C_1}{\ep^2}\EB + {C_1}'L e^{Lt}.
$$

\subsubsection {An estimate for $E_2$}

First, in the points where where $|d _\phi|<d_0$, by
\eqref{bellettini-motion-aniso}, we have $E_2=0$. Next we consider
the points where $|d _\phi|\geq d_0.$ We deduce from the
definition of $\Lambda _0$ in \eqref{lambda0-aniso} that
$$
\begin{array}{ll}
0\leq 2a(x,\n \dphi (x,t))&\leq (\Lambda ^0)^2|\n
\dphi (x,t)|^2\vsp\\
&\leq  (\Lambda ^0)^2 \Vert \n \dphi \Vert _\infty ^2:=D<\infty.
\end{array}
$$
Applying Lemma \ref{est-phi-aniso} yields
$$
|E_2|\leq \di{\frac{C}{\ep^2}}(1+D)e^{-\lambda|d _\phi+\ep p|/
\ep}\leq \di{\frac{C}{\ep^2}}(1+D)e^{-\lambda(d_0 / \ep-|p|)}.
$$
We remark that $0<K-1 \leq p \leq e^{LT} +K$, and suppose from now
that the following assumption holds:
\begin{equation}\label{cond3-aniso}
e^{LT}+K \leq \frac{d_0}{2\ep_0}.
\end{equation}
Then $\di{\frac{d_0}{\ep}}-|p|\geq \di{\frac{d_0}{2\ep}}$ so that,
defining $C':=C(1+D)$,
$$
|E_2|\leq \di{\frac{C'}{\ep^2}}e^{-\lambda d_0 / (2\ep)}\leq
C_2:=\di{\frac{16C'}{(e\lambda d_0)^2}}.
$$

\subsubsection {An estimate for $E_3$}
We set
$$
\mathcal G(x,t)=(d _\phi) _t(x,t)-\surm \div \big[m(x)a_p(x,\n
\dphi (x,t))\big].
$$
We recall that $\dphi \in C^{3+\vartheta,(3+\vartheta)/2}$ in a
neighborhood $\mathcal{V}$ of $\Gamma$, say
$$
\mathcal{V}=\{(x,t) \in \Q,\; |d _\phi(x,t)|<d_0\}.
$$
Combining the fact that
$$
2a(x,\n \dphi (x,t))=1 \quad\ \text{in} \; \ \mathcal V,
$$
with the definition of $\Lambda ^0$ in \eqref{lambda0-aniso}, we
see that
\begin{equation}\label{awayfromzero-aniso}
|\n \dphi|\geq {\frac{1}{\Lambda^0}}\quad\ \text{in} \; \ \mathcal
V.
\end{equation}
We also recall that $(x,p)\mapsto a(x,p)$ is of class
$C^{3+\vartheta}_{loc}$ on $\ombar \times \R ^N \setminus \{0\}$.
Since $|\n \dphi|$ is bounded away from zero, it follows that
$x\mapsto \mathcal G (x,t)$ is Lipschitz continuous on
$\mathcal{V}$. By equation \eqref{FBP-aniso}, we have that
$$
\mathcal G (x,t)=0 \quad \textrm{ on } \Gamma _t=\{x \in \om, d
_\phi(x,t)=0\},
$$
and it follows from the mean value theorem applied separately on
both sides of $\Gamma_t$ that there exists a constant $N_1$ such
that
\begin{equation}\label{presgamma-aniso}
|\mathcal G (x,t)|\leq N_1|d _\phi(x,t)| \quad \textrm{ for all
}(x,t) \in \mathcal V.
\end{equation}
Next, using Lemma \ref{aniso-lapla-aniso}, we remark that
$\mathcal G$ is bounded on $\ombar\times[0,T]\backslash
\mathcal{V}$ so that there exists a constant $N_2$ such that
\begin{equation}\label{loingamma-aniso}
\sup_{\ombar\times[0,T]\backslash \mathcal{V}}|\mathcal
G(x,t)|\leq N_2.
\end{equation}
By the inequalities \eqref{presgamma-aniso} and
\eqref{loingamma-aniso}, we deduce that
$$
|\mathcal G (x,t)| \leq N |\dphi (x,t)| \quad\ \text{in}\; Q_T,
$$
with $N:=\max(N_1,N_2/d_0)$. Applying  Lemma \ref{est-phi-aniso}
we deduce that
$$
\begin{array}{lll}
|E_3|&\leq NC\di{\frac
{|d _\phi|}{\ep}}e^{-\lambda| d _\phi/ \ep +p|}\vsp \\
&\leq NC \max_{y \in \mathbb{R}}|y|e^{-
\lambda|y +p|}\vsp \\
&\leq NC\max (|p|,\di{\frac 1 \lambda}).
\end{array}
$$
Thus, recalling that $|p|\leq e^{Lt}+K$, we obtain
$$
|E_3|\leq C_3(e^{Lt}+K)+{C_3}', $$ where $C_3:=NC$ and
${C_3}':=NC/\lambda$.

\subsubsection {Completion of the proof}
Collecting the above estimates of $E_1$, $E_2$ and $E_3$ yields
$$
{\cal L}_0 u_\ep^+\geq \frac{C_1}{\ep^2}\EB+
(L{C_1}'-{C_3})e^{Lt}-C_4,
$$
where $ C_4:=C_2+KC_3+{C_3}'$. Now, we set
$$
L:=\frac 1 T\ln \frac {d_0}{4\ep _0}\,,
$$
which, for $\ep_ 0$ small enough, validates assumptions
\eqref{ep0M-aniso} and \eqref{cond3-aniso}. If $\ep_0$ is chosen
sufficiently small (i.e. $L$ sufficiently large), we obtain, for
all $\ep \in (0,\ep _0)$,
$$
{\cal L}_0 u_\ep^+\geq (L{C_1}'-C_3)e^{Lt}-C_4\geq  \frac 1 2
L{C_1}' -C_4\geq 0.
$$
The proof of Lemma \ref{fix-aniso} is now completed.\qed

\section{Proof of Theorem
\ref{width-aniso} and Corollary
\ref{total-aniso}}\label{s:proof-aniso}

Let $\eta \in (0,\eta _0)$ be arbitrary. Choose $\beta$ and
$\sigma$ that satisfy \eqref{beta-aniso}, \eqref{sigma-aniso} and
\begin{equation}\label{eta-aniso}
\sigma \beta \leq \frac \eta 3.
\end{equation}
By the generation of interface Theorem \ref{g-th-gen-aniso}, there
exist positive constants $\ep_0$ and $M_0$ such that
\eqref{g-part1-aniso}, \eqref{g-part2-aniso} and
\eqref{g-part3-aniso} hold with the constant $\eta$ replaced by
$\sigma \beta /2$. Since, by the hypothesis \eqref{dalltint-aniso}
and the equality \eqref{link1-aniso}, $\n u_0 (x) \cdot n _\phi
(x) \neq 0$ everywhere on the initial interface $\Gamma
_0=\{x\in\om, \; u_0(x)=a\}$ and since $\Gamma _0$ is a compact
hypersurface, we can find a positive constant $M_1$ such that
\begin{equation}\label{corres-aniso}
\begin{array}{ll}\text { if } \quad \dphi(x,0)\geq \quad M_1 \ep &
\text { then } \quad
u_0(x) \geq a +M_0 \ep,\vspace{3pt}\\
\text { if } \quad \dphi(x,0) \leq -M_1 \ep & \text { then } \quad
u_0(x) \leq a -M_0 \ep.
\end{array}
\end{equation}
Now we define functions $H^+(x), H^-(x)$ by
\[
\begin{array}{l}
H^+(x)=\left\{
\begin{array}{ll}
1+\sigma\beta/2\quad\ &\hbox{if}\ \ \dphi(x,0)> -M_1\ep\\
\;\;\;\;\;\;\sigma\beta/2\quad\ &\hbox{if}\ \ \dphi(x,0)\leq
-M_1\ep,
\end{array}\right.
\vsp\\
H^-(x)=\left\{
\begin{array}{ll}
1-\sigma\beta/2\quad\ &\hbox{if}\ \ \dphi(x,0)\geq \;M_1\ep\\
\;\;\;-\sigma\beta/2\quad\ &\hbox{if}\ \ \dphi(x,0)<  \;M_1\ep.
\end{array}\right.
\end{array}
\]
Then from the above observation we see that
\begin{equation}\label{H-u-aniso}
H^-(x) \,\leq\, u^\ep(x,\mu^{-1}  \ep^2|\ln \ep|) \,\leq\, H^+(x),
\end{equation}
for almost all $x\in\om$.

Next we fix a sufficiently large constant $K>1$ such that
\begin{equation}\label{K-aniso}
U_0(-M_1+K) \geq 1-\frac {\sigma \beta}{3} \quad \text { and }
\quad U_0(M_1-K) \leq \frac {\sigma \beta}{3}.
\end{equation}
For this $K$, we choose $\ep _0$ and $L$  as in Lemma
\ref{fix-aniso}. We claim that
\begin{equation}\label{uep-H-aniso}
u_\ep^-(x,0)\leq H^-(x),\quad\ H^+(x)\leq u_\ep^+(x,0),
\end{equation}
for all $x\in\Omega$. We only prove the former inequality, as the
proof of the latter is virtually the same. Then it amounts to
showing that
\begin{equation}\label{c3-aniso}
u_\ep^-(x,0)=U_0\big(\frac {d_0(x)}{\ep}-K\big)-\sigma (\beta+\ep
^2 L) \;\leq\; H^-(x).
\end{equation}
In the range where $\dphi(x,0) < M_1 \ep$, the second inequality
in \eqref{K-aniso} and the fact that $U_0$ is an increasing
function imply
\begin{align*}
U_0\big(\frac {\dphi(x,0)}{\ep}-K\big)-\sigma (\beta+\ep ^2 L)
&\leq U_0(M_1 -K)-\sigma \beta - \sigma \ep^2 L\vsp \\
&\leq \frac {\sigma\beta} 3 -\sigma \beta\vsp \\
&\leq H^-(x).
\end{align*}
On the other hand, in the range where $\dphi(x,0) \geq M_1 \ep$,
we have
\begin{align*}
U_0\big(\frac {\dphi(x,0)}{\ep}-K\big)-\sigma (\beta+\ep^2L) &\leq
1-\sigma
\beta\\
&\leq H^-(x).
\end{align*}
This proves \eqref{c3-aniso}, so that \eqref{uep-H-aniso} is
established.

Combining \eqref{H-u-aniso} and \eqref{uep-H-aniso}, we obtain
$$
u_\ep^-(x,0)\leq u^\ep(x,\mu ^{-1} \ep ^2|\ln \ep|) \leq
u_\ep^+(x,0),
$$
for almost all $x\in\om$. Since, by Lemma \ref{fix-aniso},
$u_\ep^-$ and $u_\ep^+$ are sub- and super-solutions of Problem
$\Pe$, the comparison principle yields
\begin{equation}\label{ok-aniso}
u_\ep^-(x,t) \leq u^\ep (x,t+t^\ep) \leq u_\ep^+(x,t),
\end{equation}
for almost all $(x,t) \in Q_T$ that satisfies $0 \leq t \leq
T-t^\ep$, where we recall that $t^\ep=\mu ^{-1} \ep ^2|\ln \ep|$.
Note that, in view of \eqref{sub2-aniso}, this is sufficient to
prove Corollary \ref{total-aniso}. Now let $C$ be a positive
constant such that
\begin{equation}\label{C}
U_0(C-e^{LT}-K) \geq 1-\frac \eta 2 \quad \text { and } \quad
U_0(-C+e^{LT}+K) \leq \frac \eta 2.
\end{equation}

One then easily checks, using \eqref{ok-aniso} and
\eqref{eta-aniso}, that, for $\ep _0$ small enough and for almost
all $(x,t)\in \Q$, we have
\begin{equation}\label{correspon-aniso}
\begin{array}{ll}\text { if } \quad \dphi(x,t) \geq \quad C \ep &
\text { then } \quad
u^\ep (x,t+t^\ep) \geq 1 -\eta,\vspace{3pt}\\
\text { if } \quad \dphi(x,t) \leq -C \ep & \text { then } \quad
u^\ep (x,t+t^\ep) \leq \eta,
\end{array}
\end{equation}
and
$$
u^\ep (x,t+t^\ep) \in [-\eta,1+\eta],
$$
which completes the proof of Theorem \ref{width-aniso}.\qed

\end{document}